\documentclass{siamart}

\makeatletter
\def\ps@pprintTitle{%
 \let\@oddhead\@empty
 \let\@evenhead\@empty
 \def\@oddfoot{}%
 \let\@evenfoot\@oddfoot}
\makeatother

\usepackage{amsmath,amscd,amsfonts,amssymb,mathtools}
\usepackage{mathabx,breqn,multirow,comment}
\usepackage{booktabs,cool}
\usepackage{parskip}
\usepackage[top=1.2in, bottom=1.2in, left=1in, right=1in]{geometry}

\DeclareMathOperator{\sign}{sign}

\newtheorem{remark}{Remark}

\let\sep=\null
\newcommand{\maxdegree}{1\sep,000\sep,000}
\newcommand{\tablesize}{138}

\providecommand{\e}[1]{\ensuremath{\times 10^{#1}}}

\title
{
An algorithm 
for the numerical evaluation of the  associated Legendre functions
that runs in time independent of degree and order
}

\author{James Bremer\thanks{Department of Mathematics, 
University of California, Davis (bremer@math.ucdavis.edu)}
}

\begin{document}

\maketitle

\begin{abstract}
We describe a method for the numerical evaluation of normalized
versions of the associated Legendre functions $P_\nu^{-\mu}$
and $Q_\nu^{-\mu}$ of degrees $0 \leq \nu \leq \maxdegree$ 
 and orders $-\nu \leq \mu \leq \nu$ on the interval $(-1,1)$.
Our algorithm, which runs in time independent
of $\nu$ and $\mu$,  is based on the fact
that while  the associated Legendre functions themselves are extremely expensive to represent
via polynomial expansions, the logarithms of certain solutions of 
the differential equation defining them are not.    We exploit
this by numerically precomputing the logarithms of carefully
chosen solutions of the associated Legendre differential equation 
and representing them  via piecewise trivariate Chebyshev expansions.  
These precomputed expansions, which allow for  the rapid evaluation
of the associated Legendre functions over a large swath of parameter domain
mentioned above, 
are supplemented with asymptotic and series expansions in order to cover it entirely.
The results of  numerical experiments demonstrating the efficacy
of our approach are presented, and 
our code for evaluating the associated Legendre functions  is publicly available.

\end{abstract}










In this article, we describe an algorithm for the numerical evaluation of 
the functions $\tilde{P}_\nu^{-\mu}$ and $\tilde{Q}_\nu^{-\mu}$ defined via
\begin{equation}
\tilde{P}_\nu^{-\mu}(t) = 
\bar{P}_\nu^{-\mu} \left( \cos(t) \right) \sqrt{\sin(t)}
\label{introduction:p}
\end{equation}
and
\begin{equation}
\tilde{Q}_\nu^{-\mu}(t) = 
\bar{Q}_\nu^{-\mu} \left( \cos(t) \right) \sqrt{\sin(t)},
\label{introduction:q}
\end{equation}
where $\bar{P}_\nu^{-\mu}$ and $\bar{Q}_\nu^{-\mu}$
are the normalized associated Legendre functions 
\begin{equation}
\bar{P}_\nu^{-\mu}(x) = \sqrt{\left(\nu+\frac{1}{2}\right)\frac{\Gamma(\nu+\mu+1)}{\Gamma(\nu-\mu+1)}} 
\ P_\nu^{-\mu}(x)
\label{introduction:normp}
\end{equation}
and
\begin{equation}
\bar{Q}_\nu^{-\mu}(x) = 
\frac{2}{\pi}
\sqrt{\left(\nu+\frac{1}{2}\right)\frac{\Gamma(\nu+\mu+1)}{\Gamma(\nu-\mu+1)}} 
\ Q_\nu^{-\mu}(x).
\label{introduction:normq}
\end{equation}
It runs in time independent of degree $\nu$ and order $\mu$, and is applicable
when $0 \leq \nu \leq \maxdegree$, $0 \leq \mu \leq \nu$
and $0 < t \leq \frac{\pi}{2}$ (in particular, $\nu$ and $\mu$ need not be integers).
It is a consequence of standard connection formulas (such as those appearing in 
Section~3.4 of \cite{HTFI}) that this suffices for the evaluation of 
\begin{equation}
\bar{P}_\nu^{\mu}(x),\  \bar{P}_\nu^{-\mu}(x),\  \bar{Q}_\nu^{\mu}(x),\  \mbox{and}\ \  \bar{Q}_\nu^{-\mu}(x)
\label{introduction:others}
\end{equation}
for any  $0 \leq \nu \leq \maxdegree$, $-\nu \leq \mu \leq \nu$
and $-1 < x < 1$.

Our choice of scaling and normalization factors  in (\ref{introduction:normp}) 
and (\ref{introduction:normq}) are fairly standard.  Among other things, they ensure that 
\begin{equation}
\int_{-1}^1 \left( \tilde{P}_\nu^{-\mu}(t) \right) ^2 \ dt = 1
\end{equation}
whenever $-\nu \leq \mu \leq \nu$  with $\nu$ and $\mu$ integers,
and that  $\tilde{P}_\nu^{-\mu}$ and $\tilde{Q}_\nu^{-\mu}$ satisfy the second order linear ordinary 
differential equation
\begin{equation}
y''(t) + \left( \lambda^2 - \eta^2 \csc^2(t) \right) y(t) = 0
\ \ \ \mbox{for all} \ \ 0 < t < \frac{\pi}{2}
\label{introduction:diffeq}
\end{equation}
with $\lambda = \nu + \frac{1}{2}$ and $\eta^2 = \mu^2-\frac{1}{4}$.
By a slight abuse of terminology, we will refer to (\ref{introduction:diffeq})
as the associated Legendre differential equation.

When $0 \leq \mu \leq \frac{1}{2}$, the coefficient of $y$
in (\ref{introduction:diffeq}) is positive on the 
interval $\left(0,\frac{\pi}{2}\right)$, whereas when $\mu > \frac{1}{2}$
it is negative on the interval
\begin{equation}
\left(0,\arcsin\left(\frac{\eta}{\lambda}\right)\right)
\label{introduction:intervalno}
\end{equation}
and positive on
\begin{equation}
\left(\arcsin\left(\frac{\eta}{\lambda}\right),\frac{\pi}{2}\right).
\label{introduction:intervalosc}
\end{equation}
It follows from these observations and 
well-known WKB estimates (see, for example, \cite{Fedoryuk})
that when $\mu > \frac{1}{2}$, the solutions of (\ref{introduction:diffeq}) behave roughly like
combinations of increasing or decreasing exponentials on (\ref{introduction:intervalno})
and are oscillatory on (\ref{introduction:intervalosc}), 
whereas when $\mu \leq \frac{1}{2}$, they are   oscillatory on all of $\left(0,\frac{\pi}{2}\right)$.
We will refer to the subset
\begin{equation}
\begin{aligned}
\mathcal{O} &= 
\left\{
\left(\nu,\mu,t\right) :
\nu \geq 0, \ 0 \leq \mu  \leq \min\left\{\nu,\frac{1}{2}\right\}\  \mbox{and}\ 0 < t \leq \frac{\pi}{2}
\right\}\bigcup\\
&
\left\{
\left(\nu,\mu,t\right) :
\nu > \frac{1}{2} , \ \ \frac{1}{2} < \mu  \leq \nu \ \ \mbox{and}\ 
 \arcsin\left(\frac{\eta}{\lambda}\right) \leq t \leq \frac{\pi}{2}
\right\}
\end{aligned}
\end{equation}
of $\mathbb{R}^3$ as the oscillatory region, and to the subset
\begin{equation}
\mathcal{N} = 
\left\{
\left(\nu,\mu,t\right) :
\nu \geq 0, \ 
\mu > \frac{1}{2} \ \ \mbox{and} \ \
0 < t < \arcsin\left(\frac{\eta}{\lambda}\right)
\right\}
\end{equation}
as the nonoscillatory region. 
When $\nu \gg \mu > \frac{1}{2}$, the
 solutions of (\ref{introduction:diffeq}) are highly oscillatory
on (\ref{introduction:intervalosc}), and when $\nu \geq \mu \gg \frac{1}{2}$,
they behave roughly like combinations of rapidly decreasing and increasing
exponentials on (\ref{introduction:intervalno}).  
Consequently,  they cannot be effectively 
represented via polynomial expansions in the variables $\nu$, $\mu$ and $t$
on either of the sets $\mathcal{O}$ or $\mathcal{N}$,
at least for large values of  the parameters.

Nonetheless, the logarithms of certain solutions of (\ref{introduction:diffeq})
can be represented efficiently via polynomial expansions
on the sets $\mathcal{N}$ and $\mathcal{O}$.
This observation  is related to the well-known 
fact that the associated Legendre differential equation admits a nonoscillatory phase function.
Many special functions of interest posses this property as well,
at least in an asymptotic sense \cite{Miller,DLMF}.  However, 
the sheer effectiveness with which nonoscillatory phase functions
can represent solutions of the general equation
\begin{equation}
y''(t) + \lambda^2 q(t) y(t) = 0 \ \ \mbox{for all}\ \ a < t <b
\label{introduction:second_order}
\end{equation}
in  which the coefficient $q$ is smooth and positive 
appears to have been overlooked.  Indeed, under mild conditions on $q$, 
it is shown in  \cite{Bremer-Rokhlin} that there exist a positive real number $\sigma$,
 a nonoscillatory function $\alpha$  and a basis
of solutions $\{u,v\}$ of (\ref{introduction:second_order}) such that
\begin{equation}
u(t) = \frac{\cos\left(\alpha(t)\right)}{\sqrt{\left|\alpha'(t)\right|}} + O\left(\exp(-\sigma\lambda)\right)
\end{equation}
and
\begin{equation}
v(t) = \frac{\sin\left(\alpha(t)\right)}{\sqrt{\left|\alpha'(t)\right|}} + O\left(\exp(-\sigma\lambda)\right).
\end{equation}
The constant $\sigma$ is a measure of the extent to which $q$ oscillates, with
larger values of $\sigma$ corresponding to greater smoothness on the part of $q$.
The function $\alpha$ is nonoscillatory in the sense that it can be represented
using various series expansions the number of terms in which do not vary
with $\lambda$.  That is, $O(\exp(-\sigma\lambda))$ accuracy is obtained using
an $O(1)$-term expansion.
The  results of \cite{Bremer-Rokhlin} are akin to standard results on
WKB approximation in that they apply to the more general case in which
$q$ varies with the parameter $\lambda$ assuming only that
$q$ satisfies certain innocuous hypotheses independent of 
$\lambda$.  An effective numerical algorithm for the computation
of nonoscillatory phase functions  for fairly general second order differential equations
is described in \cite{BremerKummer}, although we will not need it here.
We will instead use specialized formulas which apply only in the case of associated 
Legendre functions.  However, the existence of the algorithm of \cite{BremerKummer}
and results of \cite{Bremer-Rokhlin} mean that the approach of this
paper can be applied to a large class of special functions satisfying
second order differential equations.  

The algorithm of this paper operates by numerically calculating the logarithms
of certain solutions of the associated Legendre differential equation.
We represent them  via trivariate Chebyshev expansions 
---
or rather, truncated version of these expansions 
which we call ``compressed'' trivariate Chebyshev expansions ---
the coefficients of which are stored in a table on the disk.  This table
is computed only once and is loaded into memory and used
to evaluate the associated Legendre functions rapidly.
The table used in the experiments described in this paper is 
approximately  \tablesize\ MB in size. 
We supplement these precomputed expansions
with series and asymptotic expansions  in order to cover the entire
parameter domain mentioned above.  
In addition to the values of  the functions $\tilde{P}_\nu^{-\mu}$
and $\tilde{Q}_\nu^{-\mu}$, our algorithm also produces
the values of a nonoscillatory phase function for (\ref{introduction:diffeq})
and its derivative when $(\nu,\mu,t)$ is in the oscillatory region $\mathcal{O}$
and the values of the logarithms of
$\tilde{P}_\nu^{-\mu}$ and $\tilde{Q}_\nu^{-\mu}$ when   $(\nu,\mu,t)$ is in the nonoscillatory region
$\mathcal{N}$.
The phase function is useful for, among other things, calculating the roots
of the associated Legendre functions and 
applying special function transforms involving the associated Legendre functions.
Calculating the values of the logarithms 
of (\ref{introduction:normp}) and (\ref{introduction:normq})
obviates many problems which arise from numerical overflow and underflow.

There is an extensive literature on the asymptotic behavior
of the  associated Legendre functions, and it is often suggested that existing asymptotic expansions, 
particularly Liouville-Green type expansions, suffice for the numerical evaluation of them.
While many highly effective approximations are available, it appears to be quite difficult
to produce a numerical algorithm  which is efficient and accurate
over the entire range of the variables $\nu$, $\mu$ and $t$ considered here.
The  trigonometric expansions
\begin{dmath}
\tilde{P}_\nu^{-\mu}(t) 
=
\sqrt{
\frac
{2\left(\nu + \frac{1}{2}\right) \Gamma(\nu+\mu+1)\Gamma(\nu-\mu+1)}
{\pi \sin(t) }
}
\frac{1}{\Gamma\left(\nu+\frac{3}{2}\right)} \times \\
\sum_{k=0}^\infty 
\left(
(-1)^k
\frac{\left(\mu + \frac{1}{2} \right)_k \left(-\mu + \frac{1}{2} \right)_k}
{ \left(2\sin(t)\right)^k\Gamma(k+1) \left(\nu + \frac{3}{2}\right)_k} 
\sin\left(\left(\nu+k+\frac{1}{2}\right)t + 
\frac{\pi}{2} \left( k - \mu  \right) + \frac{\pi}{4} 
\right)\right)
\label{introduction:trigp}
\end{dmath}
and
\begin{dmath}
\tilde{Q}_\nu^{-\mu}(t)
=
\sqrt{
\frac
{2\left(\nu + \frac{1}{2}\right) \Gamma(\nu+\mu+1)\Gamma(\nu-\mu+1)}
{\pi \sin(t) }
}
\frac{1}{\Gamma\left(\nu+\frac{3}{2}\right)} \times \\
\sum_{k=0}^\infty 
\left(
(-1)^k
\frac{\left(\mu + \frac{1}{2} \right)_k \left(-\mu + \frac{1}{2} \right)_k}
{ \left(2\sin(t)\right)^k\Gamma(k+1) \left(\nu + \frac{3}{2}\right)_k} 
\cos\left(\left(\nu+k+\frac{1}{2}\right)t + 
\frac{\pi}{2} \left( k - \mu  \right) + \frac{\pi}{4} 
\right)\right),
\label{introduction:trigq}
\end{dmath}
which can be found in Section~3.5 of \cite{HTFI},
illustrate some of  difficulties that arise.  
In (\ref{introduction:trigp}) and (\ref{introduction:trigq}),
$\left(x\right)_k$ is the Pochhammer symbol defined via
\begin{equation}
\left(x\right)_k = \frac{\Gamma(x+k)}{\Gamma(x)}.
\end{equation}
These expansions are applicable only when $\frac{\pi}{6} < \theta < \frac{5\pi}{6}$,
and  they  require 
on the order of $|\mu|$  terms in order to achieve a fixed accuracy.
Even more seriously, when  $|\mu|$  is not small relative to $\nu$,
the terms in  (\ref{introduction:trigp}) and (\ref{introduction:trigq}) 
are of large magnitude and alternate signs, with the consequence
that the numerical evaluation 
of (\ref{introduction:trigp}) and (\ref{introduction:trigq})
generally results in catastrophic cancellation errors.

Liouville-Green methods can be used to obtain asymptotic expansions
of the associated Legendre functions
which  are uniform in the argument $t$ and apply in the event that
 $0 \leq \mu \leq \nu$. In \cite{Boyd-Dunster}, the new dependent variable $\zeta$ defined via
the implicit relations
\begin{equation}
\int_{\gamma^2}^\zeta \frac{\sqrt{\xi^2-\gamma^2}}{2\xi }\ d\xi
= - 
\int_{\sqrt{1-\gamma^2}}^x \frac{\sqrt{1-\gamma^2-t^2}}{1-t^2}\ dt
\ \ \ \mbox{if}\ \ x < \sqrt{1-\gamma^2}
\label{introduction:zeta1}
\end{equation}
and
\begin{equation}
\int_{\gamma^2}^\zeta \frac{\sqrt{\xi^2-\gamma^2}}{2\xi }\ d\xi
= - 
\int_{\sqrt{1-\gamma^2}}^x \frac{\sqrt{t^2-1+\gamma^2}}{1-t^2}\ dt
\ \ \ \mbox{if}\ \ x > \sqrt{1-\gamma^2},
\label{introduction:zeta2}
\end{equation}
where $\lambda = \nu+\frac{1}{2}$ and $\gamma = \frac{\mu}{\lambda}$,
is introduced to obtain the uniform asymptotic expansions
\begin{equation}
\begin{aligned}
\bar{P}_\nu^{-\mu}(x) \approx
\left(\frac{\zeta-\gamma^2}{1-\gamma^2-x^2}\right)^{\frac{1}{4}}
\left(
\sqrt{\lambda}
J_{\mu} \left(\lambda \sqrt{\zeta} \right)
\sum_{k=0}^\infty  \frac{F_k(\zeta)}{\lambda^{2s}} 
+
\left(\frac{\zeta}{\lambda}\right)^{\frac{1}{2}}
J_{\mu}'\left(\lambda \sqrt{\zeta}\right)
\sum_{k=0}^\infty  \frac{\tilde{F}_k(\zeta)}{\lambda^{2s}} 
\right)
\end{aligned}
\label{introduction:BoydDunsterp}
\end{equation}
as $\nu \to\infty$ and
\begin{equation}
\begin{aligned}
\bar{Q}_\nu^{-\mu}(x) \approx
-
\left(\frac{\zeta-\gamma^2}{1-\gamma^2-x^2}\right)^{\frac{1}{4}}
\left(
\sqrt{\lambda}
 Y_{\mu } \left(\lambda \sqrt{\zeta}\right)
\sum_{k=0}^\infty  \frac{F_k(\zeta)}{\lambda^{2s}}  +
\left(\frac{\zeta}{\lambda}\right)^{\frac{1}{2}}
Y_{\mu}'\left(\lambda \sqrt{\zeta}\right)
\sum_{k=0}^\infty  \frac{\tilde{F}_k(\zeta)}{\lambda^{2s}} 
\right)
\end{aligned}
\label{introduction:BoydDunsterq}
\end{equation}
as $\nu\to\infty$.
The coefficients are given by $F_0(\zeta)=1$ 
and the formulas
\begin{equation}
\tilde{F}_k(\zeta) = \left|\zeta-\alpha^2\right|^{-\frac{1}{2}}
\int_{\alpha^2}^\zeta
\left|\xi-\alpha^2\right|^{-\frac{1}{2}}
\left( \xi F_k''(\xi)  + F_k'(\xi)  - \psi(\xi)F_k(\xi)\right) \ d\xi
\label{introduction:int1}
\end{equation}
and
\begin{equation}
F_{k+1}(\zeta) = 
- \zeta \tilde{F}_k'(\zeta)
+
\int_{\alpha^2}^\zeta
\psi(\xi) \tilde{F}_{k}'(\xi)\ d\xi 
+
 \zeta \tilde{F}_k'(\alpha^2),
\label{introduction:int2}
\end{equation}
where
\begin{equation}
\psi(\zeta) = 
\frac{1}{16}
\left(\zeta -\alpha ^2\right)^{-2}
\left(\zeta + 4 \alpha ^2
+ \frac{1-x^2}{\zeta}
\left(\frac{\zeta-\alpha^2}{1-\alpha^2-x^2}\right)^3
\left(
\left(1-4\alpha^2\right)x^2+(\alpha^4-1)
\right)
\right).
\label{introduction:psi}
\end{equation}
Note that the variable $x$  in (\ref{introduction:psi}) depends on $\zeta$ through either
(\ref{introduction:zeta1}) or (\ref{introduction:zeta2});
if this is neglected and $x$ is treated as a constant in 
 (\ref{introduction:int1}) and (\ref{introduction:int2}), then these integrals diverge.
While (\ref{introduction:BoydDunsterp}) and (\ref{introduction:BoydDunsterq})
are powerful expansions, it is not clear how to accurately and 
rapidly compute the variable $\zeta$  defined via the relations
(\ref{introduction:zeta1}) and (\ref{introduction:zeta2}) given $x$.
Nor is it obvious how to evaluate the coefficients in these expansions.
Only the first few are known analytically, and the numerical
calculation of the remaining coefficients is 
complicated by the delicate cancellations of singularities which occur 
in these formulas.
Alternate Liouville-Green expansions for  (\ref{introduction:normp})
and (\ref{introduction:normq}) are given in Chapter~12 of \cite{Olver}.  
However, the coefficients in these expansions appear to be
 no easier to compute than those in   (\ref{introduction:BoydDunsterp})
and  (\ref{introduction:BoydDunsterq}), and they also involve
a variable transformation defined implicitly by a nonlinear equation. 
Moreover, they have the unfortunate property that are only applicable when $\mu$ is small relative
to $\nu$.

The literature on asymptotic expansions of associated Legendre functions
is vast, and the possibility of constructing a numerical code 
for the evaluation of (\ref{introduction:normp}) and (\ref{introduction:normq})
using such methods  cannot be ruled out.  However, the approach
offered in this paper --- namely, the precomputation of expansions
representing the logarithms of certain solutions of the associated Legendre
differential equation ---
is simple-minded and highly  effective.
Moreover, the same basic technique can be applied, with little effort, to a large class
of special functions satisfying second order differential equations.
Indeed, in \cite{BremerBessel}, similar techniques were used to numerically precompute
a table which allows for the evaluation
the Bessel functions of the first and second kinds on the interval $(0,\infty)$
in time independent of order.

The remainder of this paper is structured as follows.  In Section~\ref{section:preliminaries},
 we review certain mathematical facts and numerical procedures which are used in the
rest of this  article.   
In Section~\ref{section:solver}, we describe a numerical
method for the solution of the differential equation
(\ref{introduction:diffeq}) which runs in time independent of the parameters $\nu$ and $\mu$.   
The construction of a precomputed table of expansions
of the associated Legendre functions which makes use
of the algorithm of Section~\ref{section:solver}
is discussed in  Section~\ref{section:expansions}.
In Section~\ref{section:algorithm}, we  detail our algorithm
for the numerical calculation of the associated Legendre functions.
Section~\ref{section:experiments} describes the results of numerical
experiments carried out to verify the efficacy of our algorithm.
We close with a few brief comments in Section~\ref{section:conclusion}.

\begin{section}{Mathematical and numerical preliminaries}

\begin{subsection}{The condition number of the evaluation of a function}

The condition number of the evaluation of a differentiable
function $f:\mathbb{R} \to\mathbb{R}$
at the point $x$  is commonly defined to be 
\begin{equation}
\kappa_{f}(x) = \left| \frac{x f'(x) }{f(x)} \right|
\label{preliminaries:condition:1}
\end{equation}
(see, for instance, Section~1.6 of \cite{Higham}). 
This quantity measures the ratio of the magnitude of 
the relative change in $f(x)$ induced by a small change
in the argument $x$ to the magnitude of the relative change in $x$
in the sense that 
\begin{equation}
\left| \frac{f(x+\delta) - f(x)}{f(x)} \right|
\approx  \kappa_{f}(x)\ \left|\frac{\delta}{x}\right|
\end{equation}
for small $\delta$.  Since almost all quantities which arise in the course
of numerical calculations  are subject to perturbations
with relative magnitudes on the order of machine epsilon,
we consider 
\begin{equation}
\kappa_f(x) \epsilon_0,
\end{equation}
where $\epsilon_0$ denotes machine epsilon, to be a rough estimate of the 
relative accuracy one should expect when evaluating $f(x)$ numerically
(in fact, it tends to be a  slightly pessimistic estimate).
In the rest of this paper, we take $\epsilon_0$ to be 
\begin{equation}
\epsilon_0 = 2^{-52} \approx 2.22044604925031\e{-16}.
\end{equation}
It is immediately clear from  (\ref{preliminaries:condition:1}) that
when $f'(x_0) x_0 \neq 0$ and $f(x_0) = 0$, $\kappa_{f}(x)$ diverges to $\infty$ as  $x \to x_0$.
One consequence of this is that there is often a significant loss
of relative accuracy when  a function is evaluated near one of its roots.
For the most part, we avoid this issue by 
representing the solutions of the associated Legendre differential equation
via functions which are bounded away from $0$.

\label{preliminaries:condition}
\end{subsection}

\begin{subsection}{Trivariate Chebyshev expansions}

For each nonnegative integer $n$, the Chebyshev polynomial of degree $n$
is defined for $-1 <x < 1$ via the formula
\begin{equation}
T_n(x) = \cos\left(n \arccos(x)\right).
\label{preliminaries:chebyshev:1}
\end{equation}
The trivariate Chebyshev series of a continuous function $f:[-1,1]^3 \to \mathbb{R}$ is 
\begin{equation}
\sideset{}{'}\sum_{i=0}^\infty  \sideset{}{'}\sum_{j=0}^\infty 
\sideset{}{'}\sum_{k=0}^\infty
a_{i,j,k} T_i(x) T_j(y) T_k(z),
\label{preliminaries:chebyshev:2}
\end{equation}
where  the coefficients are defined via the formula
\begin{equation}
a_{i,j,k} = 
\frac{8}{\pi^3}\int_{-1}^1 \int_{-1}^1 \int_{-1}^1 f(x,y,z) T_i(x) T_j(y) T_k(z) 
\frac{dx}{\sqrt{1-x^2}}\frac{dy}{\sqrt{1-y^2}} \frac{dz}{\sqrt{1-z^2}}
\label{preliminaries:chebyshev:3}
\end{equation}
and the dashes next to the summation symbols indicate that the first  term in each
sum is halved.    The well-known relationship between Chebyshev and Fourier series
(see, for instance, \cite{Mason-Hanscomb}),
together with the results of 
\cite{Fefferman1} on the pointwise almost everywhere convergence of multiple Fourier series 
immediately imply that 
\begin{equation}
\lim_{N\to\infty} \sideset{}{'}\sum_{i=0}^N \sideset{}{'}\sum_{j=0}^N
\sideset{}{'}\sum_{k=0}^N
a_{i,j,k} T_i(x) T_j(y) T_k(z) = f(x,y,z)
\label{preliminaries:chebyshev:4}
\end{equation}
for almost all $(x,y,z) \in [-1,1]^3$.
As in the case of univariate Chebyshev series,
under mild smoothness conditions on $f$, the convergence of (\ref{preliminaries:chebyshev:2}) 
is uniform.  See, for instance, Theorem~5.9 in  \cite{Mason-Hanscomb}.

If  $f(x,y,z)$ is analytic on the set 
\begin{equation}
\left\{
\left(x,y,z\right) \in \mathbb{C}^3 :\
\left|x+\sqrt{x^2-1}\right| < r_1,\  \ 
\left|y+\sqrt{y^2-1}\right| < r_2 ,\  \
\left|z+\sqrt{z^2-1}\right| < r_3 
\right\},
\label{preliminaries:chebyshev:5}
\end{equation}
where $r_1, r_2, r_3 >1$, then $\left|a_{i,j,k}\right| = O \left(r_1^{-i} r_2^{-j}  r_3^{-k}\right)$
(this is according to Theorem~11 in Chapter~V of \cite{Bochner-Martin}),  with the consequence
that the limit in (\ref{preliminaries:chebyshev:4}) converges rapidly to $f$
when $f$ is analytic in a large neighborhood containing $[-1,1]^3$.

For each nonnegative integer $n$, we refer to the collection of points
\begin{equation}
\rho_{j,n} =  -\cos\left(\frac{\pi j}{n}\right), \ \ j=0,1,\ldots,n,
\label{preliminaries:chebyshev:nodes}
\end{equation}
as the $(n+1)$-point Chebyshev grid on the interval $[-1,1]$,
and we call individual elements of this set  Chebyshev nodes or points.
One discrete version of  the well-known orthogonality relation
\begin{equation}
\int_{-1}^1   \frac{T_i(x) T_j(x)}{\sqrt{1-x^2}}\ dx = 
\begin{cases}
0   & \mbox{if}\ \ i \neq j \\
\frac{\pi}{2} & \mbox{if}\ \ i = j > 0 \\
\pi & \mbox{if}\ \ i=j=0.
\end{cases}
\label{preliminaries:chebyshev:6}
\end{equation}
is
\begin{dmath}
\sideset{}{''}\sum_{l=0}^{n} T_i(\rho_{l,n}) T_j(\rho_{l,n})
=
\begin{cases}
0   & \mbox{if}\ \ 0 \leq i,j \leq  n\ \ \mbox{and}\ i \neq j \\
\frac{n}{2} & \mbox{if}\ \ 0 < i=j <n\\
n & \mbox{if}\ \ i=j=0\ \  \mbox{or}\ \ i=j=n.
\end{cases}
\label{preliminaries:chebyshev:7}
\end{dmath}
Here, the double dash next to the summation sign
indicates that the first and last term in the series are halved.
Formula~(\ref{preliminaries:chebyshev:7})
 can be found in a slightly different form in Chapter~4 of \cite{Mason-Hanscomb}.

It follows easily from (\ref{preliminaries:chebyshev:6}) and (\ref{preliminaries:chebyshev:7})
that any trivariate polynomial $f$ of degree less than or equal to $n$ can be represented in the form
\begin{equation}
f(x,y,z) = 
 \sideset{}{''}\sum_{i=0}^n \sideset{}{''}\sum_{j=0}^n \sideset{}{''}\sum_{k=0}^n
b_{i,j,k} T_i(x) T_j(y) T_k(z),
\label{preliminaries:chebyshev:8}
\end{equation}
where
\begin{equation}
b_{i,j,k} = 
\frac{8}{n^3}
\sideset{}{''} \sum_{r=0}^n  \sideset{}{''}\sum_{s=0}^n  \sideset{}{''}\sum_{t=0}^n
T_i\left(\rho_{r,n}\right) T_j\left(\rho_{s,n}\right) T_k\left(\rho_{t,n}\right) 
f\left(\rho_{r,n}, \rho_{s,n},\rho_{t,n}\right).
\label{preliminaries:chebyshev:9}
\end{equation}
%
If $f$ is not a  polynomial of degree less than or equal to $n$, then the representation
(\ref{preliminaries:chebyshev:8}) is no longer exact.  However, in this event,
there is a well-known relationship between the coefficients
 defined via (\ref{preliminaries:chebyshev:9})
and those given by (\ref{preliminaries:chebyshev:2}).  In particular,
\begin{equation}
b_{i,j,k} =
a_{ijk} + \sum_{l_1=1}^\infty \sum_{l_2=1}^\infty \sum_{l_3=1}^\infty
\left(
a_{i+2l_1 n,j+2l_2 n,k+2 l_3 n} + a_{-i+2l_1n,-j+2l_2n,-k+2l_3 n}\right)
\label{preliminaries:chebyshev:10}
\end{equation}
for all $0 \leq i,j,k \leq n$ (the one-dimensional version of this
result can be found, for instance, in \cite{Mason-Hanscomb}).
Using (\ref{preliminaries:chebyshev:10}) it is easy to show that 
there exists a constant $C$ such that
\begin{dmath}
\sup_{x \in [-1,1]}
\left|
f(x,y,z) - \sideset{}{''} \sum_{i=0}^n \sideset{}{''} \sum_{j=0}^n \sideset{}{''} \sum_{k=0}^n
b_{i,j,k} T_i(x) T_j(y) T_k(z)\right|
\leq
C \sup_{x \in [-1,1]}
\left|
f(x,y,z) - \sideset{}{''} \sum_{i=0}^n \sideset{}{''} \sum_{j=0}^n \sideset{}{''} \sum_{k=0}^n
a_{i,j,k} T_i(x) T_j(y) T_k(z)\right|.
\label{preliminaries:chebyshev:11}
\end{dmath}
It follows, of course, that the sum (\ref{preliminaries:chebyshev:9}) converges rapidly 
to $f$ when $f$ is analytic in a large neighborhood of $[-1,1]^3$.
By a slight abuse of terminology, we will refer to (\ref{preliminaries:chebyshev:8})
as the $n^{th}$ order Chebyshev expansion for the function $f$.

\label{preliminaries:chebyshev}
\end{subsection}


\begin{subsection}{Compressed trivariate Chebyshev expansions}

It often happens that  many of the coefficients in the trivariate Chebyshev expansion
(\ref{preliminaries:chebyshev:8}) of a function $f:[-1,1]^3\to\mathbb{R}$
are of negligible magnitude.  In order to reduce the cost of storing 
such expansions as well as the cost of evaluating them, we use the following
construction to reduce the number of coefficients which need to be considered.

Suppose that $\epsilon > 0$, and that 
\begin{equation}
\sideset{}{''}\sum_{i=0}^n \sideset{}{''}\sum_{j=0}^{n}  \sideset{}{''}\sum_{k=0}^{n}
b_{i,j,k} T_i(x) T_j(y) T_k(z)
\label{preliminaries:compressed:1}
\end{equation}
is the $n^{th}$ order Chebyshev expansion for $f:[-1,1]^3 \to \mathbb{R}$.
We let $M$ denote the least nonnegative integer which is less than
or equal to $n$ and such that
\begin{equation}
\left|b_{i,j,k}\right| < \epsilon
\ \ \mbox{for all} \ \ i=M+1,\ldots,n,\ \ j=0,\ldots,n\ \ \mbox{and} \ \ k=0,\ldots,n,
\end{equation}
assuming such an integer exists.  If not, then we take $M=n$.
For each $i=0,\ldots,M$, we let $m_i$ be the least nonnegative integer
less than or equal to $n$ such that
\begin{equation}
\left|b_{i,j,k}\right| < \epsilon
\ \ \mbox{for all}\ \  j=m_i+1,\ldots,n, \ \ \mbox{and}\ \ k=0,\ldots,n
\end{equation}
if such an integer exists, and we let $m_i = n$ otherwise.
Finally, for each pair $(i,j)$ such that
$0 \leq i \leq M$ and $0 \leq j \leq m_i$, we let $n_{i,j}$
be the least nonnegative integer  such that
\begin{equation}
\left|b_{i,j,k}\right| < \epsilon
\ \ \mbox{for all}\ \   k=0,\ldots,n_{i,j}.
\end{equation}
We refer to the series
\begin{equation}
\sum_{i=0}^M
\sum_{j=0}^{m_i} 
\sum_{k=0}^{n_{i,j}}\widetilde{b}_{i,j,k} T_i(x) T_j(x) T_k(x),
\label{preliminaries:chebyshev4:expansion}
\end{equation}
where $\widetilde{b}_{i,j,k}$ is defined via
\begin{equation}
\widetilde{b}_{i,j,k} = b_{i,j,k} 
\left( 1 - \frac{1}{2}  \delta_{i,n}\right)\left( 1  - \frac{1}{2} \delta_{i,0}\right)  
\left( 1  - \frac{1}{2} \delta_{j,0}\right)
\left( 1  - \frac{1}{2} \delta_{j,n}\right)
\left( 1  - \frac{1}{2} \delta_{k,0}\right)
\left( 1  - \frac{1}{2} \delta_{k,n}\right),
\end{equation}
as the $\epsilon$-compressed $n^{th}$ order Chebyshev expansion of $f$.

Obviously, the results discussed in this and the preceding section
 can be modified in 
a straightforward fashion so as to apply to functions given on an arbitrary 
rectangular prism $[a,b] \times [c,d] \times [e,f]$.
%

\label{preliminaries:compressed}
\end{subsection}


\begin{subsection}{Series expansions of the associated Legendre functions and connection formulas}

When $\nu \geq 0$, $-\nu \leq \mu \leq \nu$ and $-1 \leq x < 1$, the associated Legendre function of the 
first kind of degree $\nu$ and order $-\mu$ is given by 
\begin{equation}
P_\nu^{-\mu}(x) = \left(\frac{1-x}{1+x}\right)^{\mu/2}
\sum_{n=0}^\infty  (-1)^n \frac{\Gamma(\nu + n + 1)}{\Gamma(\nu - n +1 )}
\frac{\left(\frac{1}{2} - \frac{x}{2}  \right)^n}
{\Gamma(n+1)\Gamma(n+\mu+1)}.
\label{preliminaries:series:1}
\end{equation}
Here, we have adopted the convention that 
\begin{equation}
\frac{1}{\Gamma(k)} = 0
\end{equation}
whenever $k$ is a negative integer.   The trigonometric form
\begin{equation}
P_\nu^{-\mu}(\cos(t)) = 
\left( \tan\left(\frac{t}{2}\right) \right)^{\mu}
\sum_{n=0}^\infty \left(-1\right)^n \frac{\Gamma(\nu + n + 1)}{\Gamma(\nu - n +1 )}
\frac{\left(\sin\left(\frac{t}{2}\right)\right)^{2n}}
{\Gamma(n+1)\Gamma(n+\mu+1)}
\label{preliminaries:series:2}
\end{equation}
 of (\ref{preliminaries:series:1}) is obtained by letting $x=\cos(t)$ and
making use of elementary identities.   When 
the parameters $\nu$ and $\mu$ are of small magnitude,
the coefficients in  (\ref{preliminaries:series:2}) decay rapidly as $n$ increases,
with the consequence that only a small number of
terms of (\ref{preliminaries:series:2})
are   required to accurately evaluate $P_\nu^{-\mu}$.  
Likewise, even when the parameters are of large magnitude
 the coefficients  in this expansion decay rapidly with $n$ 
 if $t$ is sufficiently small, so that
(\ref{preliminaries:series:2}) is efficient in this regime as well.
For extremely large values of $\nu$, we found
 numerical roundoff error to be a problem
in the evaluation of (\ref{preliminaries:series:2}).  For this reason,
we only use this series expansion  in the event that $\nu$ 
is less than $10\sep,000$.

One potential difficulty with the use of (\ref{preliminaries:series:2}) as a numerical tool,
however, is that underflow can occur when the parameters are large
and $t$ is small.  To obviate such problems, 
we use a truncation of the formula
\begin{dmath}
\log\left(\tilde{P}_\nu^{-\mu}(t)\right) = 
\log\left(\frac{\Gamma(\nu-\mu+1)}{\Gamma(\nu+\mu+1)}\right)
+\log\left(\nu+\frac{1}{2}\right)
+ \frac{1}{2}\log\left(\sin(t)\right)
+\mu \log\left(\tan\left(\frac{t}{2}\right) \right)
- \log\left(\frac{1}{\Gamma\left(\mu+1\right)}\right)
+
\log\left(
\sum_{n=0}^\infty \left(-1\right)^n \frac{\Gamma(\nu + n + 1)}{\Gamma(\nu - n +1 )}
\frac{\Gamma(\mu+1)}{\Gamma(\mu+n+1)}
\frac{\left(\sin\left(\frac{t}{2}\right)\right)^{2n}}
{\Gamma(n+1)}
\right),
\label{preliminaries:series:2.5}
\end{dmath}
which is easily obtained from  (\ref{introduction:p}) and (\ref{preliminaries:series:2}),
to evaluate the logarithm of the associated Legendre function of the first kind
in this regime.  We note that $\tilde{P}_\nu^{-\mu}(t)$ is necessarily positive
when $t$ is sufficiently small, so that this logarithm is sensible.
When $\mu$ is equal to a negative integer, say $-m$, 
the first $m$ terms of the sum in (\ref{preliminaries:series:2})
are $0$ and we use a version of (\ref{preliminaries:series:2.5}) which is modified
accordingly.    

\vskip 1em
\begin{remark}
 The naive evaluation of the first term in 
(\ref{preliminaries:series:2.5}) can lead to numerical cancellation
when $\nu$ is large and $\mu$ is small relative to $\nu$.  In this event, 
we use the first sixteen terms of the 
asymptotic approximation
\begin{dmath*}
\log\left(\frac{\Gamma(x-y)}{\Gamma(x+y)}\right) \approx 2\left(\mu+1\right)\log\left(\nu\right) 
+ \log \left( 1
-\frac{y}{x}
-\frac{y \left(2 y^2-3   y+1\right)}{6 x^2} 
+\frac{y^2 \left(2 y^2-3 y+1\right)}{6 x^3}+
\frac{y \left(20 y^5-96y^4+155 y^3-90 y^2+5 y+6\right)}{360 x^4}
+\cdots\right)
\end{dmath*}
in order to evaluate it.
\label{remark1}
\end{remark}

For $\nu \geq 0$, $-1 \leq x < 1$ and $-\nu \leq \mu \leq \nu$ not an integer,
the associated Legendre function of the second kind of degree $\nu$
and order $-\mu$ is given by 
\begin{equation}
Q_\nu^{-\mu}(x) = \frac{\pi}{2}
\left(
 \frac{\Gamma(\nu-\mu+1)}{\Gamma(\nu+\mu+1)} \sec(\mu\pi) P_\nu^{\mu}(x)- \cot(\mu \pi)P_\nu^{-\mu}(x)
\right).
\label{preliminaries:series:3}
\end{equation}
The normalized versions of the associated Legendre
functions satisfy the somewhat simpler relation
\begin{equation}
\tilde{Q}_\nu^{-\mu}(t) = \sec(\mu\pi)
\tilde{P}_\nu^{\mu}(t)-  \cot(\mu \pi)\tilde{P}_\nu^{-\mu}(t),
\label{preliminaries:series:4}
\end{equation}
which is an immediate consequence of (\ref{preliminaries:series:3}), 
(\ref{introduction:p}) and (\ref{introduction:q}).
Similarly to the case of (\ref{preliminaries:series:2}), 
the use of  (\ref{preliminaries:series:4}) can lead to numerical overflow
when $t$ is small.
Accordingly, we generally compute the logarithm of $\tilde{Q}_\nu^{-\mu}$
via the less delicate formula
\begin{dmath}
\log\left(\tilde{Q}_\nu^{-\mu}(t)\right) = 
\log\left(\tilde{P}_\nu^{-\mu}(t) \right)
+ \log\left(
\sec(\mu\pi)
-  \cot(\mu \pi) \sign\left(\tilde{P}_\nu^{\mu}(t)\right)
\\ \exp\left(\log\left(\left|\tilde{P}_\nu^{\mu}(t)\right|\right) - 
\log\left(\tilde{P}_\nu^{-\mu}(t)\right)\right)
\right)
\label{preliminaries:series:6}
\end{dmath}
in this regime.  We note that
for sufficiently small $t$,
the function $\tilde{Q}_\nu^{-\mu}(t)$  is positive.

When $\mu$ is an integer, (\ref{preliminaries:series:4}) 
and (\ref{preliminaries:series:6}) lose their meanings.
Various series expansions for $Q_\nu^{-m}(x)$ with  $m$ a positive integer 
  can be obtained (see, for instance, Section 3.6 of \cite{HTFI}),
but they are somewhat cumbersome and  do not address a second problem
with the use of  (\ref{preliminaries:series:4})
 as a numerical method for the evaluation of  $Q_\nu^{-\mu}$.  
Namely, that  when $\mu$ is close to, but does not coincide with, an integer,
the evaluation (\ref{preliminaries:series:4}) 
 results in severe loss of precision due to numerical cancellation.
However, since $\tilde{Q}_\nu^{-\mu}$  is an analytic function of the parameter
$\mu$, it  can be efficiently interpolated in the
$\mu$ variable.    For instance, when $\mu$ is close to, or coincides with,
an integer $m$, the value of   $\tilde{Q}_\nu^{-\mu}(t)$ can be calculated
by first evaluating
\begin{equation}
 \tilde{Q}_\nu^{\xi_1}(t), \ldots,  \tilde{Q}_\nu^{\xi_{2n+1}}(t)
\end{equation}
with $\xi_1,\ldots,\xi_{2n}$ the nodes of the $(2n)$-point Chebyshev grid on the interval
$[m-\epsilon],[m+\epsilon]$, and then using Chebyshev interpolation
to calculate $\tilde{Q}_\nu^{-\mu}(t)$. 
  Here, $\epsilon$ is an appropriate chosen positive
real number  and $n$ is a positive integer.  
 An even number of nodes is chosen in order to ensure that none coincide with the integer $m$.
In the code used in this paper,
we apply this procedure when $\mu$ is within a distance of $0.001$ of an
integer, and we take $n=6$ and $\epsilon=0.1$.
Of course, the same approach can be used to evaluate $\log\left(\tilde{Q}_\nu^{-\mu}(t)\right)$.

Indeed, many other connection formulas for the associated Legendre functions can be handled in a similar 
fashion, such as the identity
\begin{equation}
\tilde{Q}_\nu^{-\mu}(t) = \sec(\mu \pi) \tilde{Q}_\nu^\mu(t) + 
\tan(\mu\pi) \tilde{P}_\nu^{-\mu}(t),
\label{preliminaries:series:7}
\end{equation}
which follows easily from a formula found in 
in Section~3.4 of \cite{HTFI}.
On the other hand, the connection formulas
\begin{equation}
\tilde{P}_\nu^{-\mu}(\pi-t) = 
\cos(\pi(\nu-\mu)) \tilde{P}_\nu^{-\mu}(t)
- \sin(\pi(\nu-\mu)) \tilde{Q}_\nu^{-\mu}(t)
\end{equation}
and
\begin{equation}
\tilde{Q}_\nu^{-\mu}(\pi-t) = 
-\cos(\pi(\nu-\mu)) \tilde{Q}_\nu^{-\mu}(t)
- \sin(\pi(\nu-\mu)) \tilde{P}_\nu^{-\mu}(t),
\end{equation}
which also appear (in a slightly different form) in Section~3.4 of \cite{HTFI},
are immune from such problems.

\label{preliminaries:series}
\end{subsection}


\begin{subsection}{Macdonald's asymptotic expansions}

In \cite{Macdonald}, an asymptotic formula for $P_\nu^{-\mu}(\cos(t))$ which is accurate
when $\nu$ is large, $0 \leq \mu \leq \nu$ and $t$ is small
is derived by replacing the ratio of Gamma functions
\begin{equation}
\frac{\Gamma(\nu+n+1)}{\Gamma(\nu-n+1)}
\label{preliminaries:macdonald:1}
\end{equation}
appearing in (\ref{preliminaries:series:2})
with a finite truncation of the series expansion
\begin{equation}
\frac{\Gamma(\nu+n+1)}{\Gamma(\nu-n+1)}
= 
\lambda^{2n} 
- G_1 \lambda^{2n-2} 
+ G_2 \lambda^{2n-4}  
+ G_3 \lambda^{2n-6} + \cdots,
\label{preliminaries:macdonald:2}
\end{equation}
where $\lambda = \nu + \frac{1}{2}$.
The first three terms in the  asymptotic expansion of $P_\nu^{-\mu}$ obtained 
in this fashion
are
\begin{dmath}
P_\nu^{-\mu}(\cos(t)) \approx \left( \lambda \cos\left(\frac{t}{2}\right) \right)^{-\mu}
\cdot \left( J_\mu(\eta)  + \sin^2 \left(\frac{t}{2}\right) H_1 
+ \sin^4\left(\frac{t}{2}\right) H_2 + \sin^6\left(\frac{t}{2}\right) H_3 \right),
\label{preliminaries:macdonald:pexpansion}
\end{dmath}
where 
%
\begin{dmath*}
H_1 = \frac{\eta}{6}J_{\mu+3}(\eta)  -  J_{\mu+2}(\eta)   +  \frac{1}{2\eta} J_{\mu+1}(\eta),
\end{dmath*}
\begin{dmath*}
H_2 =\frac{\eta^2}{72} J_{\mu+6}(\eta) - \frac{11 \eta}{30} J_{\mu+5}(\eta) 
+  \frac{31}{12} J_{\mu+4}(\eta) -\frac{29}{6 \eta}  J_{\mu+3}(\eta)   + \frac{9}{8\eta^2}
J_{\mu+2}(\eta),
\end{dmath*}
and
\begin{dmath*}
H_3 = \frac{75}{16 \eta^3} J_{\mu+2}(\eta)
-\frac{751}{24 \eta^2} J_{\mu+3}(\eta)
+\frac{1381}{48 \eta} J_{\mu+4}(\eta)
-\frac{1513}{180} J_{\mu+5}(\eta)
+\frac{4943 \eta}{5040} J_{\mu+6}(\eta)
-\frac{17 \eta^2}{360} J_{\mu+7}(\eta)
+\frac{\eta^3}{1296} J_{\mu+8}(\eta).
\end{dmath*}

The first few terms of the analogous expansion
of the associated Legendre function of the second kind, which
is applicable when $0 \leq \mu \leq \nu$,  are
\begin{dmath}
Q_\nu^{\mu}(\cos(t)) \approx 
-\frac{\pi}{2}
\left( \lambda \cos\left(\frac{t}{2}\right) \right)^{\mu}
\cdot \left( Y_{-\mu}(\eta)  + \sin^2 \left(\frac{t}{2}\right) N_1
+ \sin^4\left(\frac{t}{2}\right) N_2 + \sin^6\left(\frac{t}{2}\right) N_3 \right),
\label{preliminaries:macdonald:qexpansion}
\end{dmath}
where  $Y_\zeta$ denotes the Bessel function of the second
kind of order $\zeta$, 
\begin{dmath*}
N_1 = \frac{\eta}{6}Y_{-\mu+3}(\eta)  -  Y_{-\mu+2}(\eta)   +  \frac{1}{2\eta} Y_{-\mu+1}(\eta),
\end{dmath*}
\begin{dmath*}
N_2 =\frac{\eta^2}{72} Y_{-\mu+6}(\eta) - \frac{11 \eta}{30} Y_{-\mu+5}(\eta) 
+  \frac{31}{12} Y_{-\mu+4}(\eta) -\frac{29}{6 \eta}  Y_{-\mu+3}(\eta)   + \frac{9}{8\eta^2}
Y_{-\mu+2}(\eta),
\end{dmath*}
and
\begin{dmath*}
N_3 = \frac{75}{16 \eta^3} Y_{-\mu+2}(\eta)
-\frac{751}{24 \eta^2} Y_{-\mu+3}(\eta)
+\frac{1381}{48 \eta} Y_{-\mu+4}(\eta)
-\frac{1513}{180} Y_{-\mu+5}(\eta)
+\frac{4943 \eta}{5040} Y_{-\mu+6}(\eta)
-\frac{17 \eta^2}{360} Y_{-\mu+7}(\eta)
+\frac{\eta^3}{1296} Y_{-\mu+8}(\eta).
\end{dmath*}

To reduce the potential for numerical underflow
in the evaluation of   (\ref{preliminaries:macdonald:pexpansion}),
we evaluate $\log\left(\tilde{P}_\nu^{-\mu}(t)\right)$ using the following formula instead:
\begin{dmath}
\log\left(\tilde{P}_\nu^{-\mu}(t)\right)
\approx 
\log\left(\frac{\Gamma(\nu-\mu+1)}{\Gamma(\nu+\mu+1)}\right)
+\log\left(\nu+\frac{1}{2}\right)
+ \frac{1}{2}\log\left(\sin(t)\right)
- \mu \log \left(\lambda \cos\left(\frac{t}{2}\right)\right)
+
\log\left(J_\mu(\eta)\right)
+
\log\left(
1 + \sin^2\left(\frac{t}{2}\right) \tilde{H}_1
+ \sin^4\left(\frac{t}{2}\right) \tilde{H}_2
+\sin^6\left(\frac{t}{2}\right) \tilde{H}_3
\right),
\label{preliminaries:macdonald:pexpansion_log}
\end{dmath}
where
\begin{dmath*}
\tilde{H}_1 = 
\frac{\eta}{6}
\exp\left(\log\left(J_{\mu+3}(\eta)\right)-\log\left(J_\mu(\eta)\right)\right)
  -  
\exp\left(\log\left(J_{\mu+2}(\eta)\right)-\log\left(J_\mu(\eta)\right)\right)
  +  
\frac{1}{2\eta} \exp\left(\log\left(J_{\mu+1}(\eta)\right)-\log\left(J_\mu(\eta)\right)\right),
\end{dmath*}
\begin{dmath*}
\tilde{H}_2 =
\frac{\eta^2}{72} 
\exp\left(\log\left(J_{\mu+6}(\eta)\right)-\log\left(J_\mu(\eta)\right)\right)
- \frac{11 \eta}{30} 
\exp\left(\log\left(J_{\mu+5}(\eta)\right)-\log\left(J_\mu(\eta)\right)\right)
+  \frac{31}{12} 
\exp\left(\log\left(J_{\mu+4}(\eta)\right)-\log\left(J_\mu(\eta)\right)\right)
-\frac{29}{6 \eta}  
\exp\left(\log\left(J_{\mu+3}(\eta)\right)-\log\left(J_\mu(\eta)\right)\right)
+ \frac{9}{8\eta^2}
\exp\left(\log\left(J_{\mu+2}(\eta)\right)-\log\left(J_\mu(\eta)\right)\right),
\end{dmath*}
and
\begin{dmath*}
\tilde{H}_3 = \frac{75}{16 \eta^3} 
\exp\left(\log\left(J_{\mu+2}(\eta)\right)-\log\left(J_\mu(\eta)\right)\right)
-\frac{751}{24 \eta^2} 
\exp\left(\log\left(J_{\mu+3}(\eta)\right)-\log\left(J_\mu(\eta)\right)\right)
+\frac{1381}{48 \eta} 
\exp\left(\log\left(J_{\mu+4}(\eta)\right)-\log\left(J_\mu(\eta)\right)\right)
-\frac{1513}{180} 
\exp\left(\log\left(J_{\mu+5}(\eta)\right)-\log\left(J_\mu(\eta)\right)\right)
+\frac{4943 \eta}{5040} 
\exp\left(\log\left(J_{\mu+6}(\eta)\right)-\log\left(J_\mu(\eta)\right)\right)
-\frac{17 \eta^2}{360} 
\exp\left(\log\left(J_{\mu+7}(\eta)\right)-\log\left(J_\mu(\eta)\right)\right)
+\frac{\eta^3}{1296} 
\exp\left(\log\left(J_{\mu+8}(\eta)\right)-\log\left(J_\mu(\eta)\right)\right).
\end{dmath*}
As discussed in  Remark~\ref{remark1}, some care must be taken in evaluating the 
first term in (\ref{preliminaries:macdonald:pexpansion_log}).  
We use an analogous form of (\ref{preliminaries:macdonald:qexpansion}) 
in order to evaluate $\log\left(\tilde{Q}_\nu^{\mu}(t)\right)$.  
The logarithms of the Bessel functions appearing in these formulas
are calculated  via the algorithm of \cite{BremerBessel}.

We are unaware of a simple analog of (\ref{preliminaries:macdonald:qexpansion})
for the function $\tilde{Q}_\nu^{-\mu}$.  When we say that
we evaluate $\log\left(\tilde{Q}_\nu^{-\mu}(t)\right)$ via Macdonald's expansions,
we mean that we combine the appropriate formulas for 
$\log\left(\tilde{P}_\nu^{-\mu}(t)\right)$
and $\log\left(\tilde{Q}_\nu^{\mu}(t)\right)$
with the connection formula
(\ref{preliminaries:series:7}) in order to calculate
$\log\left(\tilde{Q}_\nu^{-\mu}(t)\right)$.

\label{preliminaries:macdonald}
\end{subsection}


\begin{subsection}{Riccati's equation, Kummer's equation and phase functions}

In this section, we suppose that $q$ is a smooth, real-valued function defined on an open
interval $I \subset \mathbb{R}$.
In the event that $q$ is strictly negative on $I$, 
two linearly independent solutions of  the second order differential equation
\begin{equation}
y''(t) + q(t) y(t) = 0 \ \ \mbox{for all}\ \ t \in I
\label{preliminaries:riccati:diffeq}
\end{equation}
both of which are positive on $I$ can be found.   This follows easily
from standard proofs of Picard's theorem 
on the existence and uniqueness of solutions of ordinary differential
equations  (see, for instance, Section~2.3 of \cite{Hille}).
Any positive solution $y$ of (\ref{preliminaries:riccati:diffeq})
can be represented in the form $y=\exp(r(t))$ with $r$ real-valued, and a straightforward
computation shows that $r$ must satisfy
\begin{equation}
r''(t)  + (r'(t))^2 + q(t) = 0
\ \ \mbox{for all}\ \ t \in I.
\label{preliminaries:riccatieq}
\end{equation}
Equation~(\ref{preliminaries:riccatieq}) is known as Riccati's equation;
a detailed discussion of it can be found in \cite{Hille}, among
many other sources.

When $q$ is positive on $I$, 
the solutions of 
(\ref{preliminaries:riccati:diffeq}) oscillate and their logarithms
are complex-valued.  In this case, it is convenient to 
represent the solutions of  (\ref{preliminaries:riccati:diffeq})
via a phase function, which
is  nothing more than the imaginary part of the logarithm
of one of its solutions.
More precisely, we say that a smooth function $\alpha$ defined on $I$
is a phase function  for the second order differential equation
(\ref{preliminaries:riccati:diffeq})
provided $\alpha'$ does not vanish on $I$ and the pair
\begin{equation}
u(t) = \frac{\cos\left(\alpha(t)\right)}{\sqrt{\left|\alpha'(t)\right|}}
\label{preliminaries:riccati:u}
\end{equation}
and
\begin{equation}
v(t) = \frac{\sin\left(\alpha(t)\right)}{\sqrt{\left|\alpha'(t)\right|}}
\label{preliminaries:riccati:v}
\end{equation}
form a basis in the space of solutions of (\ref{preliminaries:riccati:diffeq}).    
We note that the definition of phase function 
does not require that $q$ be positive, although phase functions
are most useful on intervals where this is the case.
Proofs of the following elementary
results  regarding phase functions can be found in
\cite{Heitman-Bremer-Rokhlin-Vioreanu} and \cite{BremerBessel}.

\vskip 1em
\begin{theorem}
Suppose that $I$ is an open interval in $\mathbb{R}$, and that
$q$ is a smooth, real-valued function defined on $I$.  Suppose also
that $\alpha$ is a smooth, real-valued function defined on $I$ whose
first derivative does not vanish there.  Then $\alpha$
is a phase function for the second order differential
equation (\ref{preliminaries:riccati:diffeq}) if and only if
its derivative $\alpha'$
 satisfies the second order nonlinear differential equation
\begin{equation}
 q(t) 
- (\alpha'(t))^2
- \frac{1}{2}\left(\frac{\alpha'''(t)}{\alpha'(t)}\right)
+ \frac{3}{4}
\left(\frac{\alpha''(t)}{\alpha'(t)}\right)^2 = 0
\ \ \mbox{for all}\  \ t\in I.
\label{preliminaries:riccati:kummer}
\end{equation}
\label{theorem1}
\end{theorem}

\vskip 1em
\begin{theorem}
Suppose that $u,v$ is a pair of smooth, real-valued solutions
of (\ref{preliminaries:riccati:diffeq}) whose (necessarily constant)
Wronskian $W$ is nonzero.  Then  there is a phase function $\alpha$
for (\ref{preliminaries:riccati:diffeq}) such that
\begin{equation}
u(t) = \sqrt{W} \frac{\cos(\alpha(t)) }{\sqrt{|\alpha'(t)|}}
\label{preliminaries:riccati:u2}
\end{equation}
and
\begin{equation}
v(t) =  \sqrt{W}  \frac{\sin(\alpha(t)) }{\sqrt{|\alpha'(t)|}}.
\label{preliminaries:riccati:v2}
\end{equation}
Moreover, the derivative of $\alpha$ is given by
\begin{equation}
\alpha'(t) = \frac{W}{(u(t))^2 + (v(t))^2}
\ \ \mbox{for all} \ \ t \in I,
\end{equation}
and $\alpha$ is unique up to addition by an integer
multiple of $2\pi$.   That is, $\tilde{\alpha}$ is a phase
function for (\ref{preliminaries:riccati:diffeq}) 
such that (\ref{preliminaries:riccati:u2}) and (\ref{preliminaries:riccati:v2})
hold  if and only there exists an integer $L$ such that
\begin{equation}
\tilde{\alpha}(t) = \alpha(t) + 2\pi L
\ \ \mbox{for all} \ \ t \in I.
\end{equation}
\label{theorem2}
\end{theorem}

We will refer to (\ref{preliminaries:riccati:kummer}) as Kummer's equation,
after E.~E.~Kummer who studied it in \cite{Kummer}.  

\label{preliminaries:riccati}
\end{subsection}

\begin{subsection}{A nonoscillatory phase function for the associated Legendre
differential equation}

From Theorem~\ref{theorem2}, we see that there is a phase
function $\alpha_{\nu,\mu}$ for (\ref{introduction:diffeq}) such that
\begin{equation}
\tilde{P}_\nu^{-\mu}(t) = 
\sqrt{\frac{2 \left(\nu+\frac{1}{2}\right)}{\pi}}
\frac{\cos(\alpha_{\nu,\mu}(t))}{\sqrt{\alpha_{\nu,\mu}'(t)}}
\label{preliminaries:u}
\end{equation}
and
\begin{equation}
\tilde{Q}_\nu^{-\mu}(t) = 
\sqrt{\frac{2 \left(\nu+\frac{1}{2}\right)}{\pi}}
\frac{\sin(\alpha_{\nu,\mu}(t))}{\sqrt{\alpha_{\nu,\mu}'(t)}},
\label{preliminaries:v}
\end{equation}
and whose derivative is given by
\begin{equation}
\alpha_{\nu,\mu}'(t) 
=  
\frac{2}{\pi} \left(\nu+\frac{1}{2}\right)
\frac{1}{\left(\tilde{P}_\nu^{-\mu}(t)\right)^2 + \left(\tilde{Q}_\nu^{-\mu}(t)\right)^2 }.
\label{preliminaries:alphap}
\end{equation}
We have made use of  the fact (which can be found in a slightly
different form in Section~3.4 of \cite{HTFI})
that the Wronskian of  the pair $\tilde{P}_\nu^{-\mu}, \tilde{Q}_\nu^{-\mu}$ is 
$\frac{2}{\pi}\left(\nu+\frac{1}{2}\right)$.
It has long been known that the function (\ref{preliminaries:alphap})
is nonoscillatory.   Indeed, it is immediate from (\ref{introduction:BoydDunsterp}) and 
(\ref{introduction:BoydDunsterq}) that
\begin{equation}
\begin{aligned}
\alpha_{\nu,\mu}'(t) 
&\sim
\frac{2}{\pi}
\left(\frac{\zeta-\gamma^2}{1-\gamma^2-x^2}\right)^{-\frac{1}{2}}
\frac{1}{
\left( J_{\mu } \left(\lambda \sqrt{\zeta}\right)\right)^2 + 
\left( Y_{\mu } \left(\lambda \sqrt{\zeta}\right)\right)^2 
}\ \ \mbox{as}\ \ \nu\to\infty,
\end{aligned}
\label{preliminaries:nonoscillatory:asymp}
\end{equation}
where $\zeta$ is the variable defined implicitly by (\ref{introduction:zeta1})
and (\ref{introduction:zeta2}), $\lambda=\nu+\frac{1}{2}$
and $\gamma = \frac{\mu}{\lambda}$.
A cursory inspection of Nicholson's integral formula
\begin{equation}
J_\mu^2(z) + Y_\mu^2(z) = \frac{8}{\pi^2}
\int_0^\infty K_0(2z\sinh(t)) \cosh(2\mu t)\ dt,
\label{introduction:nicholson}
\end{equation}
a derivation of which can be found in Section~13.73 of \cite{Watson},
reveals that the function
$\left( J_{\mu }(x) \right)^2 + \left( Y_{\mu }(x)\right)^2$ is nonoscillatory.
%
We note that this property of $\alpha_{\nu,\mu}$ is highly unusual.
Figure~\ref{figure1} compares $\alpha_{\nu,\mu}'$
with the derivative of typical phase functions
for (\ref{introduction:diffeq}), which oscillate
on some portion of $\left(0,\frac{\pi}{2}\right)$.

\begin{figure}[b!!]
\begin{center}
\includegraphics[width=.46\textwidth]{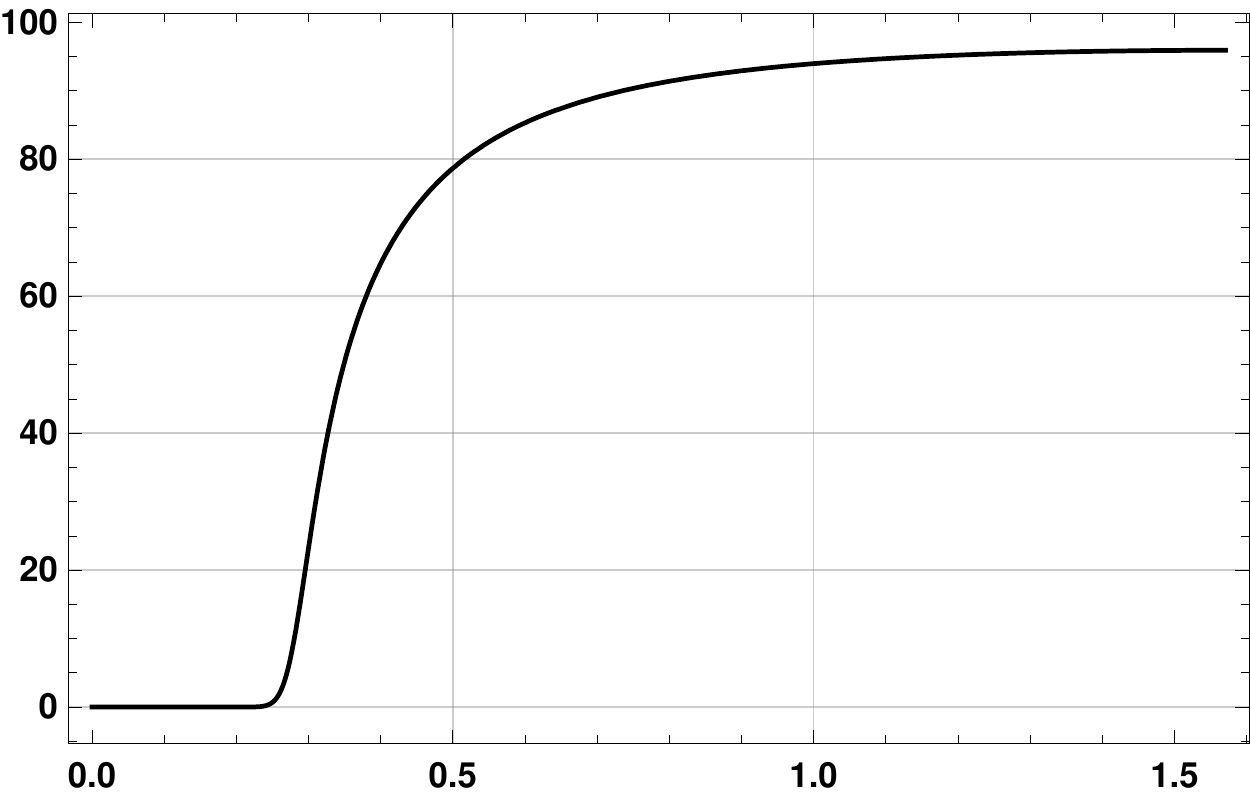}
\hfil
\includegraphics[width=.46\textwidth]{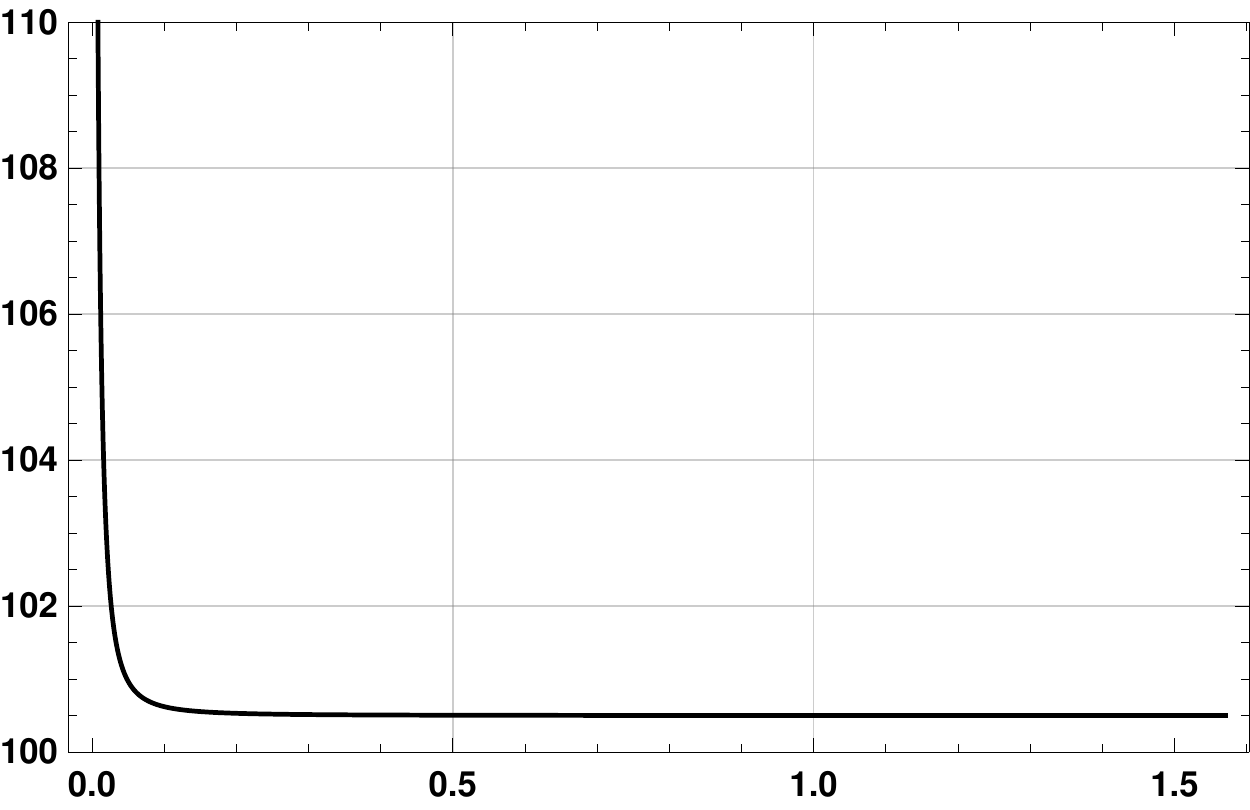}

\vskip 1.2em

\includegraphics[width=.46\textwidth]{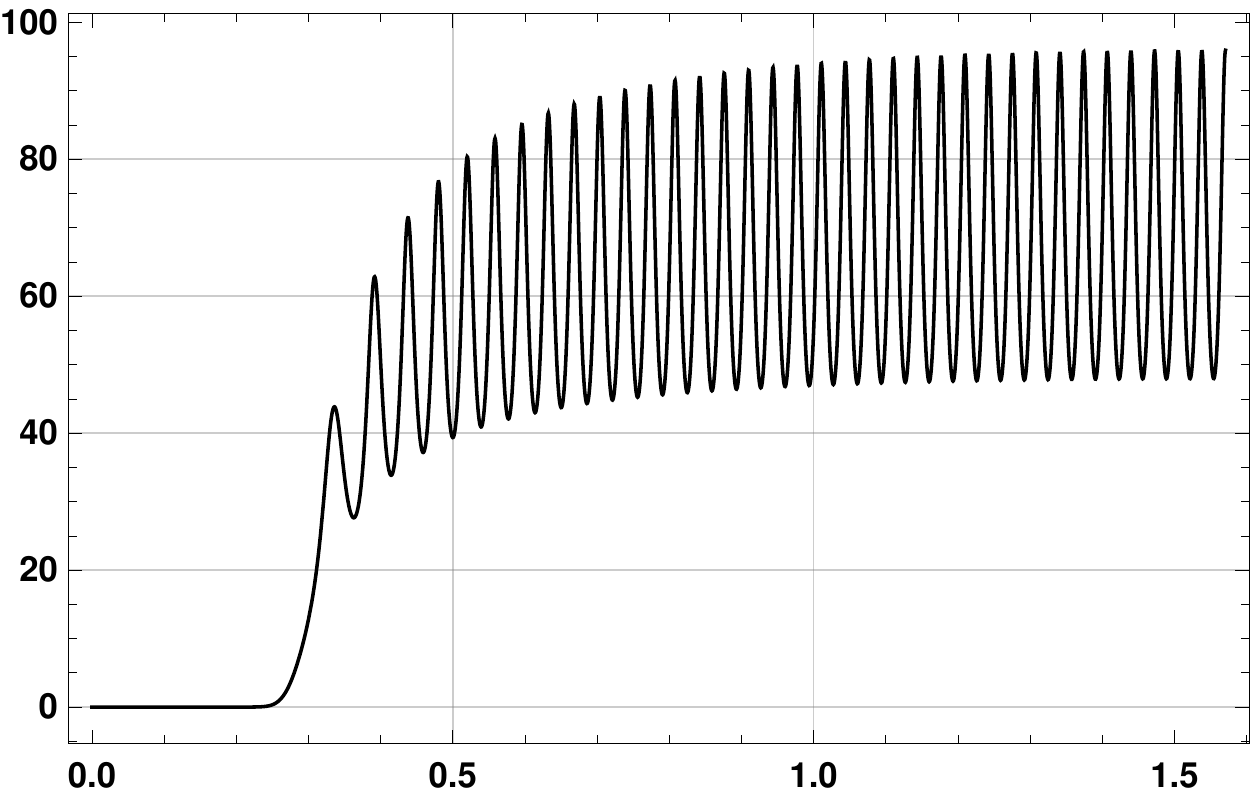}
\hfil
\includegraphics[width=.46\textwidth]{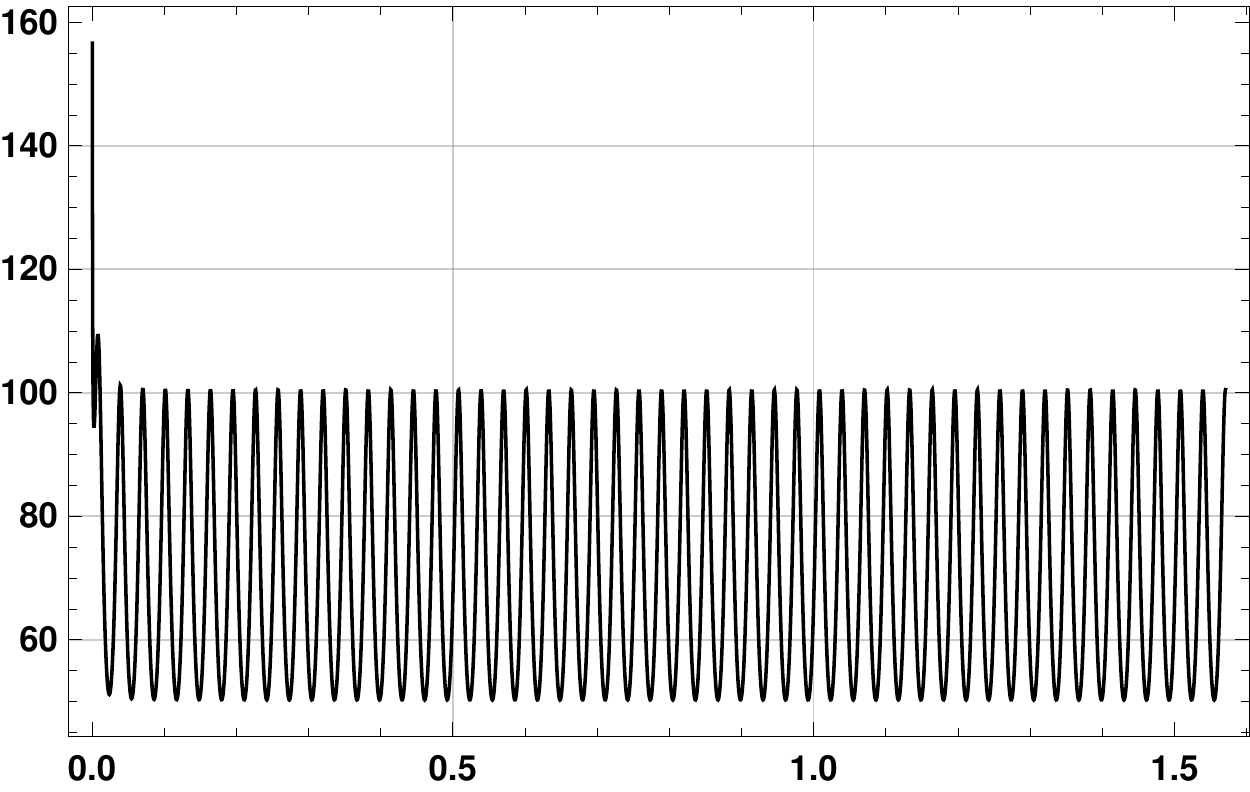}

\end{center}
\caption{
Top left: a plot of the function
 $\alpha'_{\nu,\mu}$ defined via (\ref{preliminaries:alphap}) when $\nu=100$
and $\mu=30$.  
Top right:  a plot of $\alpha'_{\nu,\mu}$ when
$\nu = 100$ and $\mu = 0$.
Bottom left: a plot of the derivative of a typical oscillatory phase
function for (\ref{introduction:diffeq}) 
when $\nu = 100$ and $\mu=30$.
Bottom right: a plot of the derivative of a  typical oscillatory phase
function for (\ref{introduction:diffeq}) 
when $\nu = 100$ and $\mu=0$.
}
\label{figure1}
\end{figure}


It follows from (\ref{preliminaries:alphap})
that  there exists a constant $C$ such that
\begin{equation}
\alpha_{\nu,\mu}(t) = C + \int_{\frac{\pi}{2}}^t \alpha_{\nu,\mu}'(s)\ ds.
\label{preliminaries:alpha0}
\end{equation}
In fact, using the formulas
\begin{equation}
\tilde{P}_\nu^{-\mu}\left(\frac{\pi}{2}\right) =  
\frac{2^{-\mu}}{\sqrt{\pi}}
\sqrt{\left(\nu+\frac{1}{2}\right)\frac{\Gamma(\nu+\mu+1)}{\Gamma(\nu-\mu+1)}} 
\frac{  \cos \left(\frac{1}{2} \pi  (\nu-\mu )\right) \Gamma
   \left(\frac{1}{2} (-\mu +\nu +1)\right)}
{ \Gamma \left(\frac{1}{2} (\mu +\nu
   +2)\right)}
\label{preliminaries:nonoscillatory:p0}
\end{equation}
and
\begin{equation}
\tilde{Q}_\nu^{-\mu}\left(\frac{\pi}{2}\right) = 
-
\frac{2^{-\mu}}{\sqrt{\pi}}
\sqrt{\left(\nu+\frac{1}{2}\right)\frac{\Gamma(\nu+\mu+1)}{\Gamma(\nu-\mu+1)}} 
 \frac{ \sin \left(\frac{1}{2} \pi  (\nu-\mu)\right) \Gamma
   \left(\frac{1}{2} (-\mu +\nu +1)\right)}{\Gamma \left(\frac{1}{2} (\mu +\nu
   +2)\right)},
\label{preliminaries:nonoscillatory:q0}
\end{equation}
which appear in a slightly different form in Section~3.4 of \cite{HTFI},
we see that (\ref{preliminaries:u}) and (\ref{preliminaries:v})
hold so long as the constant $C$ in (\ref{preliminaries:alpha0})
is taken to be 
\begin{equation}
C = \frac{\pi}{2} \left(\nu-\mu\right) + 2\pi L
\end{equation}
with $L$ an integer.   In the remainder of this paper,
we let $\alpha_{\nu,\mu}$ denote the phase function
defined via the formula
\begin{equation}
\alpha_{\nu,\mu}(t) = 2\pi + \frac{\pi}{2} \left(\nu-\mu\right) + \int_{\frac{\pi}{2}}^t \alpha_{\nu,\mu}'(s)\ ds.
\label{preliminaries:alpha}
\end{equation}
We set $L=1$ in order to  ensure that $\alpha_{\nu,\mu}$ is bounded
away from $0$ on the interval $\left(0,\frac{\pi}{2}\right)$.  In this way,
we avoid certain difficulties which arise because the condition
number of evaluation of a function is generally infinite
near one of its roots (as per the discussion in
Section~\ref{preliminaries:condition}).
%
By inserting  (\ref{preliminaries:nonoscillatory:p0}) and
 (\ref{preliminaries:nonoscillatory:q0}) 
into (\ref{preliminaries:alphap}), we see that
\begin{equation}
\alpha_{\nu,\mu}'\left(\frac{\pi}{2}\right) = 
\frac{2 \Gamma\left( \frac{1}{2} \left( \nu - \mu + 2 \right)\right)
\Gamma\left( \frac{1}{2} \left( \nu + \mu + 2 \right)\right)}
{\Gamma\left( \frac{1}{2} \left( \nu - \mu + 1 \right)\right)
\Gamma\left( \frac{1}{2} \left( \nu + \mu + 1 \right)\right)}.
\label{preliminaries:alphapright}
\end{equation}
Expressions for the values of the derivatives of the functions
$\tilde{P}_\nu^{-\mu}$ and $\tilde{Q}_\nu^{-\mu}$
at the point $\frac{\pi}{2}$ can be easily derived from
formulas appearing in
Section~3.4 of \cite{HTFI}.  A tedious computation
which makes use of them 
in addition to (\ref{preliminaries:nonoscillatory:p0}) and
 (\ref{preliminaries:nonoscillatory:q0}) 
shows that
\begin{equation}
\alpha_{\nu,\mu}\left(\frac{\pi}{2}\right) = 0.
\label{preliminaries:alphappright}
\end{equation}
%
%
\end{subsection}

\begin{subsection}{An adaptive discretization procedure}

We now briefly describe a fairly standard procedure for
 adaptively discretizing a smooth function $f:[a,b] \to \mathbb{R}$.
It takes as input a desired precision $\epsilon >0$, a positive integer $n$
and a subroutine for evaluating $f$.
The goal of this procedure is to construct a partition
\begin{equation}
a = \gamma_0 < \gamma_1 < \cdots < \gamma_m = b
\end{equation}
of $[a,b]$ such that  the $n^{th}$ order Chebyshev expansion of $f$ on each of the subintervals
$[\gamma_j,\gamma_{j+1}]$ of $[a,b]$ approximates $f$ with accuracy $\epsilon$.  That is,
for each $j=0,\ldots,m-1$ we aim to achieve
\begin{equation}
\sup_{x \in [\gamma_j,\gamma_{j+1}]}
\left| 
f(x) - 
\sideset{}{''}\sum_{i=0}^n b_{i,j} T_{i} \left( \frac{2}{\gamma_{j+1}-\gamma_j } x + 
\frac{\gamma_{j+1}+\gamma_j}{\gamma_j-\gamma_{j+1}} \right)
\right| < \epsilon,
\label{preliminaries:adaptive:1}
\end{equation}
where $b_{0,j},b_{1,j}\ldots,b_{n,j}$ are the coefficients in the $n^{th}$ order Chebyshev expansion
of $f$ on the interval $\left[\gamma_j,\gamma_{j+1}\right]$.  
These coefficients are defined by the formula
\begin{equation}
b_{i,j} = \frac{2}{n} \sideset{}{''}\sum_{l=0}^n T_i\left(\rho_{l,n}\right) 
f\left( \frac{\gamma_{j}-\gamma_{j+1}}{2}\cos\left(\frac{\pi l}{n} \right) + \frac{\gamma_{j+1}+\gamma_j}{2}\right).
\end{equation}

During the procedure,  two lists of subintervals are maintained: a list
of subintervals which are to be processed and a list of output subintervals.
Initially, the list of subintervals to be processed consists of $[a,b]$ 
and the list of output subintervals is empty.
The procedure terminates when the list of subintervals to be processed is empty
or when the number of subintervals in this list  exceeds a present
limit (we usually take this limit to be $300$).
In the latter case, the procedure is deemed to have failed.
As long as the list of subintervals to process is nonempty and its length
does not exceed the preset maximum, the algorithm
proceeds by removing a subinterval $\left[\eta_1,\eta_2\right]$
from that list  and performing the following operations:
\begin{enumerate}

\item
Compute the coefficients $b_0,\ldots,b_n$ 
in the $n^{th}$ order Chebyshev expansion
of the restriction of $f$ to the interval $\left[\eta_1,\eta_2\right]$.

\vskip 1em
\item
Compute the quantity 
\begin{equation}
\Delta = 
\frac{\max\left\{
\left|b_{\frac{n}{2}+1} \right|, 
\left|b_{\frac{n}{2}+2} \right|, \ldots
\left|b_n \right|
\right\}}
{\max\left\{
\left|b_0 \right|, 
\left|b_1 \right|, \ldots
\left|b_n \right|
\right\}}.
\end{equation}

\vskip 1em
\item
If $\Delta < \epsilon$ then 
 the subinterval $\left[\eta_1,\eta_2\right]$ is added to the list
of output subintervals.

\vskip 1em
\item 

If $\Delta \geq \epsilon$, then the subintervals
\begin{equation}
\left[\eta_1,\frac{\eta_1+\eta_2}{2}\right]
\ \ \mbox{and}\ \ 
\left[\frac{\eta_1+\eta_2}{2}, \eta_2 \right]
\end{equation}
are added to the list of  subintervals to be processed.

\end{enumerate}

This algorithm is  heuristic in the sense that there is no guarantee
that (\ref{preliminaries:adaptive:1}) will be achieved, but 
similar adaptive discretization procedures are widely used
with great success. 

There is one common circumstance which 
leads to the failure of this procedure.  The quantity $\Delta$ is an attempt to estimate
the relative accuracy with which the Chebyshev expansion
of $f$ on the interval $\left[\eta_1,\eta_2\right]$ approximates
$f$.  In cases in which the condition number of the evaluation of
$f$ is larger than $\epsilon$ on some part of $[a,b]$,
the procedure will generally fail or an excessive number of subintervals
will be generated.   Particular care needs to be taken when $f$ has a zero
in $[a,b]$.
In most cases, for $x$ near  a zero of $f$, the condition number of evaluation
of $f(x)$ (as defined in Section~\ref{preliminaries:condition})
is large.   In this article, we avoid such difficulties  by only
applying this procedure to 
 functions which are bounded away from $0$.

\label{preliminaries:adaptive}
\end{subsection}

\label{section:preliminaries}
\end{section}

\begin{section}{A method for the rapid numerical solution of the associated Legendre differential equation}

In this section, we describe an algorithm for the numerical solution
of the associated Legendre differential equation which runs in time
independent of $\nu$ and $\mu$.    It is a crucial component
of the scheme of the following section for the construction
of a table which allows for the rapid numerical evaluation
of the associated Legendre functions.  

The algorithm makes use of a solver for 
nonlinear second order ordinary differential equations 
of the form
\begin{equation}
y''(t) = f(t, y(t),y'(t)) \ \ \mbox{for all} \ \ a < t < b
\label{solver:equation}
\end{equation}
which is described in detail in  Section~4 of \cite{BremerBessel}.  
That solver is designed to be extremely robust, but not necessarily
highly efficient.    
It  takes as input a subroutine for evaluating
the function $f$ and its derivatives with respect to $t$, $y$ 
and $y'$, a positive integer $k$, a precision $\epsilon > 0$ for the calculations,
and either initial or terminal conditions for the desired
solution $y$.  It returns a collection of subintervals
\begin{equation}
\left[\gamma_1,\gamma_2\right], \ldots, \left[\gamma_{m-1},\gamma_m\right]
\label{solver:subintervals}
\end{equation}
and the values of the functions $y$, $y'$ and $y''$ at the $(k+1)$-point
Chebyshev grid on each of the subintervals (\ref{solver:subintervals}).
In particular, the functions $y$, $y'$ and $y''$ are represented
via piecewise $k^{th}$ order Chebyshev expansions.  
Given this data, the value of any one of these functions at any point on the 
interval $(a,b)$ can be computed using Chebyshev
interpolation (see, for instance, \cite{Trefethen} for a thorough discussion
of such techniques).  The collection
of subintervals is determined adaptively in the course of solving
(\ref{solver:equation}) using an approach  which
attempts to achieve relative accuracy 
in the expansions of $y$, $y'$ and $y''$ on the order of the specified precision $\epsilon$.  
The algorithm is heuristic and offers no
accuracy guarantees, but similar approaches are commonly 
used with great success.  There is one situation in which this solver
 is likely to fail.  
When the condition number
of evaluation of the solution of (\ref{solver:subintervals}) is large,
it is not possible to represent it with high relative accuracy
using Chebyshev expansions.
In this event, the solver tends
to produce an excessive number of subintervals or fail altogether.
Since the condition number of evaluation of a function is generally
large near one of its roots, we only apply this solver in cases in which
the solution is bounded away from $0$.

Our algorithm for the numerical solution of (\ref{introduction:diffeq})
takes as input real numbers $\nu$ and $\mu$ such that
$0 \leq \mu \leq \nu$, a desired precision $\epsilon > 0$,
and a positive integer $k$ specifying the order of the Chebyshev
expansions to use.  
It proceeds in three stages.  

\vskip 1em
{\it Stage one: computation of the nonoscillatory phase function $\alpha_{\nu,\mu}$}

In this stage, we calculate the values of the nonoscillatory phase
function (\ref{preliminaries:alpha}) on the interval
\begin{equation}
\left(t^*_{\nu,\mu},
\frac{\pi}{2} \right),
\label{solver:oscint}
\end{equation}
where
\begin{equation}
t^*_{\nu,\mu}  =
\arcsin\left(\frac{\sqrt{\mu^2-\frac{1}{4}}}{\nu+\frac{1}{2}}\right)
\label{solver:tstar1}
\end{equation}
is the turning point of (\ref{introduction:diffeq}) 
if $\mu \geq 1 $ and
\begin{equation}
t^*_{\nu,\mu}  = \frac{1}{\nu^{\frac{3}{2}}}
\label{solver:tstar2}
\end{equation}
if  $0 \leq \mu < 1$.
The rationale for using (\ref{solver:tstar2}) as the left
endpoint for the interval on which the phase function
is calculated when $0 < \mu < 1$ is to 
 avoid a discontinuity in $t_{\nu,\mu}$
when $\mu$ crosses the threshold $\mu = \frac{1}{2}$.

We first construct
$\alpha_{\nu,\mu}'$ by solving a terminal value problem for Kummer's equation
(\ref{preliminaries:riccati:kummer}) using the solver of Section~4 of \cite{BremerBessel}.
The values of $\alpha_{\nu,\mu}'$ and $\alpha_{\nu,\mu}''$ at the
point $\frac{\pi}{2}$  are obtained using (\ref{preliminaries:alphapright})
and (\ref{preliminaries:alphappright}).  
The required precision for these computations is taken to be  $\epsilon$. 
Next,
$\alpha_{\nu,\mu}$ is constructed through Formula~(\ref{preliminaries:alpha}).
Since $\alpha_{\nu,\mu}$ is represented via its values at the $(k+1)$-point 
Chebyshev nodes on a collection of intervals, it is easy to evaluate
the required integral via spectral integration.

Upon the completion of this stage, the values of $\alpha_{\nu,\mu}$
and its first two derivatives are known at the nodes of the
$(k+1)$-point Chebyshev points on each interval in a  collection 
of subintervals of (\ref{solver:oscint}).  Using standard
Chebyshev interpolation methods, the values of these functions
can be calculated in a stable fashion anywhere on the interval (\ref{solver:oscint}).

\vskip 1em
{\it Stage two: computation of $\log\left(\tilde{Q}_\nu^{-\mu}(t)\right)+\nu$}

In the event that $\mu  \geq 1$, we calculate the function
$\log\left(\tilde{Q}_\nu^{-\mu}(t)\right)+\nu$ on the interval
\begin{equation}
\left(
\frac{t^*_{\nu,\mu}}{100},t^*_{\nu,\mu}\right)
\label{solver:intnonsc}
\end{equation}
by solving a terminal boundary value problem for Riccati's equation
(\ref{preliminaries:riccatieq}) using the solver described in
Section~4 of \cite{BremerBessel}.
In fact, we solve the terminal boundary value problem on the slightly
larger interval
\begin{equation}
\left(
\frac{t_{\nu,\mu}^*}{100}, t_0
\right),
\end{equation}
where $t_0$ is the solution of the nonlinear equation
\begin{equation}
\alpha_{\nu,\mu}\left(t_0\right) = \frac{\pi}{4} + 2\pi.
\label{solver:neq}
\end{equation}
The functions $\alpha_{\nu,\mu}$ and its derivative having been
calculated in the preceding stage, there is no difficulty
in using Newton's method to solve (\ref{solver:neq}).
From (\ref{preliminaries:v}) and (\ref{solver:neq}), we see that
\begin{equation}
\tilde{Q}_\nu^{-\mu}\left(t_0\right) = 
\sqrt{\frac{ \left(\nu+\frac{1}{2}\right)}{\pi  \alpha_{\nu,\mu}'(t)}}
\label{solver:3000}
\end{equation}
and
\begin{equation}
\frac{d\tilde{Q}_\nu^{-\mu}}{dt} \left(t_0\right) = 
\sqrt{\frac{\nu+\frac{1}{2}}{\pi}}
\left(
\sqrt{\alpha'\left(t_0\right)}
- \frac{\alpha''\left(t_0\right)}{2 \left(\alpha'\left(t_0\right)\right)^{\frac{3}{2}}}
\right).
\label{solver:3001}
\end{equation}
The rationale for introducing $t_0$ is to ensure that the terminal
value of $\tilde{Q}_\nu^{-\mu}$ and its derivative used 
in the solution of Riccati's equation are computed
accurately.
The condition number of evaluation of the function $\tilde{Q}_\nu^{-\mu}$
is large when the parameters $\nu$ and $\mu$ are of large
magnitude, with the consequence that its numerical evaluation 
will generally result  in a loss of precision in this event.
In the case of (\ref{preliminaries:v}), the evaluation
of a trigonometric functions at a large argument is the specific mechanism
by which this loss of precision takes place.
By evaluating  $\tilde{Q}_\nu^{-\mu}$ at a point $t_0$ at which
the value of the phase function is known, however,
we avoid this loss of precision entirely.  This can be seen
from (\ref{solver:3000}) and (\ref{solver:3001}).
They involve only the evaluation of $\alpha_{\nu,\mu}'$
and $\alpha_{\nu,\mu}''$,  the condition number of evaluation
of which is  small independent of $\nu$ and $\mu$.

We construct 
$\log\left(\tilde{Q}_\nu^{-\mu}(t)\right)+\nu$ in lieu of
$\log\left(\tilde{Q}_\nu^{-\mu}(t)\right)$  because the former
is bounded away from $0$ on the interval  
(\ref{solver:intnonsc}) while the latter is not.  
Upon the completion of this stage, the values of 
$\log\left(\tilde{Q}_\nu^{-\mu}(t)\right)+\nu$
 at the nodes of the
$(k+1)$-point Chebyshev points on each of a collection
of subintervals which cover (\ref{solver:oscint}) are known.  Using standard
Chebyshev interpolation methods, the values of this function
can be calculated in a stable fashion anywhere on the interval 
(\ref{solver:intnonsc}).

\vskip 1em
{\it Stage three: computation of $\log\left(\tilde{P}_\nu^{-\mu}(t)\right)-\nu$}

Assuming that $\mu \geq 1$, we now compute 
$\log\left(\tilde{P}_\nu^{-\mu}(t)\right)-\nu$ on the interval
(\ref{solver:intnonsc}).  Proceeding here as we did in
 the calculation of  $\log\left(\tilde{Q}_\nu^{-\mu}(t)\right)+\nu$
would be problematic.  Unlike 
$\log\left(\tilde{Q}_\nu^{-\mu}(t)\right)+\nu$, which
increases rapidly as $t$ goes to $0$ from the right,
$\log\left(\tilde{P}_\nu^{-\mu}(t)\right)-\nu$
converges to $0$ rapidly as $t$ goes to $0$ from the right.
Consequently, it is recessive when solving Riccati's equation
in the backward direction and attempts to 
 approximate it numerically by solving a terminal 
value problem for Riccati's  equation lead to excessively
large errors.

Instead, we solve an initial value problem for Riccati's equation
on the interval (\ref{solver:intnonsc}) in order
to calculate $\log\left(\tilde{P}_\nu^{-\mu}(t)\right)-\nu$.
This is numerically viable  since it is a dominant solution
of  Riccati's equation when solving in the forward direction.
When $\nu \leq 10\sep,000$, we use a truncation of 
the series expansion (\ref{preliminaries:series:2.5}) in order to 
generate the necessary initial values.  
For $\nu > 10\sep,000$, we calculate initial values via
(\ref{preliminaries:macdonald:pexpansion_log})  instead
since (\ref{preliminaries:series:2.5}) can lead 
 to numerical roundoff errors when $\nu$ is large.

As before, the rationale for computing 
$\log\left(\tilde{P}_\nu^{-\mu}(t)\right)-\nu$ in lieu of
$\log\left(\tilde{P}_\nu^{-\mu}(t)\right)$  is that the former
is bounded away from $0$ on the interval  
(\ref{solver:intnonsc}) while the latter is not.  
Upon the completion of this stage, the values of 
$\log\left(\tilde{P}_\nu^{-\mu}(t)\right)-\nu$
 at the nodes of the
$(k+1)$-point Chebyshev points on each interval in a  collection 
of subintervals of (\ref{solver:oscint}) are known.  Using standard
Chebyshev interpolation methods, the values of this function
can be calculated in a stable fashion anywhere on the interval 
(\ref{solver:intnonsc}).

\vskip 1em
\begin{remark}
Although the algorithm described in this section is highly specialized
to the case of the associated Legendre differential equation, it can, in fact, be modified
so as to apply to a large class of second order differential equations
of the form
\begin{equation}
y''(t) + q(t) y(t) = 0 \ \ \mbox{for all} \  \ a < t <b .
\label{phase:1000}
\end{equation}
Suppose, for instance, that  $q$ is smooth on $[a,b]$, has  a zero at $t_{00} \in (a,b)$, 
is negative on  $(a,t_{00})$  and is positive on  $(t_{00},b)$. 
The procedure of the first stage for constructing
a nonoscillatory phase function on $(t_{00},b)$ relies on an asymptotic expansion
which allows for the evaluation
of a nonoscillatory phase function at the point $b$.
In the absence of such an approximation, the algorithm of
\cite{BremerKummer} can be used instead.  That algorithm also
proceeds by solving Kummer's equation,
but it incorporates a mechanism for numerically calculating
the appropriate initial values of a nonoscillatory phase function and its
derivatives.

The procedure of the second stage does not rely on any asymptotic or series
expansions of associated Legendre functions, only on the values of the phase function
computed in the first phase.  Consequently, it does not need to be
modified in order to obtain a solution of Riccati's equation
 which is increasing as $t \to 0^+$.

In the third stage, a series or asymptotic expansion
is used to compute the values of $\tilde{P}_\nu^{-\mu}$
and its derivative at a point near $0$.  In the event that such an approximation is 
not available, a solution of Riccati's equation which is increasing
as $t \to t_{00}$ from the left can be obtained by solving an initial
value problem with arbitrary initial conditions
and then scaling the result in order to make it consistent with 
the desired solution of (\ref{phase:1000}).
This procedure is analogous to that used in order to obtain
a recessive solution of a linear recurrence relation by running
the recurrence relation backwards (see, for instance, Section 3.6 of \cite{DLMF}).
\end{remark}

\label{section:solver}
\end{section}

\begin{section}{The numerical construction of the precomputed table}

In this section, we describe the procedure used to construct the precomputed
table which allows for the
rapid numerical evaluation of the associated Legendre functions  $\tilde{P}_\nu^{-\mu}$ and 
$\tilde{Q}_\nu^{-\mu}$ for a large range of $\nu$, $\mu$ and $t$.
This table stores the coefficients in the compressed piecewise trivariate Chebyshev
expansions 
 of eights pair of functions.

A first pair of functions $A_1$, $B_1$ allows for the 
evaluation of the phase function $\alpha_{\nu,\mu}$ and its derivative
on the subset
\begin{equation}
\mathcal{O}_1 = \left\{
(\nu,\mu,t) : 
10 \leq \nu  \leq \maxdegree,\ \  1 \leq \mu \leq \nu
\ \mbox{and}\ \ t^*_{\nu,\mu} \leq t \leq \frac{\pi}{2}
\right\}
\end{equation}
of the oscillatory region $\mathcal{O}$.  Here, $t_{\nu,\mu}^*$
is as in (\ref{solver:tstar1}) and (\ref{solver:tstar2}).
A second pair  $A_2, B_2$ allows for the evaluation
of the phase function $\alpha_{\nu,\mu}$ and its derivative
on the subset
\begin{equation}
\mathcal{O}_2 = \left\{
(\nu,\mu,t) : 
10 \leq \nu  \leq \maxdegree,\ \  0 \leq \mu < 1
\ \mbox{and}\ \ t^*_{\nu,\mu} \leq t \leq \frac{\pi}{2}
\right\}
\end{equation}
of the oscillatory region $\mathcal{O}$.  The functions
$A_3$ and $B_3$ allow for the evaluation of $\alpha_{\nu,\mu}$
and $\alpha_{\nu,\mu}'$ on 
\begin{equation}
\mathcal{O}_3 = \left\{
(\nu,\mu,t) : 
2 \leq \nu  \leq 10,\ \  1 \leq \mu \leq \nu
\ \mbox{and}\ \ t^*_{\nu,\mu} \leq t \leq \frac{\pi}{2}
\right\},
\end{equation}
and a fourth pair $A_4, B_4$ 
allows for the evaluation of the phase function $\alpha_{\nu,\mu}$
and its derivative on
\begin{equation}
\mathcal{O}_4 = 
\left\{
(\nu,\mu,t) : 
2 \leq \nu  \leq 10,\  0 \leq \mu < 1
\ \mbox{and}\ \ t^*_{\nu,\mu} \leq t \leq \frac{\pi}{2}
\right\},
\end{equation}
We divide the range of the parameter $\nu$ because
 it is more efficient to represent $\alpha_{\nu,\mu}$ via polynomial expansions
in $\frac{1}{\nu}$ when $\nu$ is large, and
via expansions in $\nu$ when $\nu$ is small.

A fifth set of functions $C_1$ and $D_1$ allows for the evaluation
of the functions 
\begin{equation}
\log\left(\tilde{P}_{\nu}^{-\mu}(t)\right)-\nu
\ \ \mbox{and}\ \ 
\log\left(\tilde{Q}_{\nu}^{-\mu}(t)\right)+\nu
\label{expansions:10001}
\end{equation}
on the subset
\begin{equation}
\mathcal{N}_1 = 
\left\{
(\nu,\mu,t) : 
10 \leq \nu  \leq \maxdegree,\  \  1 < \mu \leq \nu
\ \mbox{and}\ \ 0 \leq t < t^*_{\nu,\mu} 
\right\}
\end{equation}
of the nonoscillatory region $\mathcal{N}$.  
A sixth pair of functions $C_2$ and $D_2$ allows for the evaluation
of the functions  (\ref{expansions:10001})
on 
\begin{equation}
\mathcal{N}_2 = 
\left\{
(\nu,\mu,t) : 
10 \leq \nu  \leq \maxdegree,\  \  0 \leq \mu < 1
\ \mbox{and}\ \ 0 \leq t < t^*_{\nu,\mu} 
\right\}.
\end{equation}
The seventh pair of functions $C_3$, $D_3$ allows for the evaluation
of (\ref{expansions:10001})
 on 
\begin{equation}
\mathcal{N}_3 = 
\left\{
(\nu,\mu,t) : 
2 \leq \nu  \leq 10,\  \  1 \leq \mu \leq \nu
\ \mbox{and}\ \ 0 \leq t < t^*_{\nu,\mu} 
\right\}.
\end{equation}
The eighth and final pair of functions $C_4$, $D_4$ allows
for the evaluation of the functions
(\ref{expansions:10001}) on
\begin{equation}
\mathcal{N}_4 = 
\left\{
(\nu,\mu,t) : 
2 \leq \nu  \leq 10,\  \  0 \leq \mu < 1
\ \mbox{and}\ \ 0 \leq t < t^*_{\nu,\mu} 
\right\}.
\end{equation}
We construct expansions of the functions
(\ref{expansions:10001})
rather than expansions of
 $\log\left(\tilde{P}_{\nu}^{-\mu}(t)\right)$ 
and $\log\left(\tilde{Q}_{\nu}^{-\mu}(t)\right)$ because
the former are bounded away from $0$ on the sets in
which we consider them while the latter are not.
This ensures that their condition number of evaluation
is not large because of the presence of roots.

These computations were conducted in IEEE quadruple precision
arithmetic in order to ensure high accuracy.  The resulting
table, which consists of the coefficients in the expansions
of the functions $A_1, \ldots, A_4$, $B_1,\ldots,B_4$, $C_1,\ldots,C_4$,
$D_1,\ldots,D_4$ is roughly \tablesize\ MB in size.  
The precomputed table allows for the evaluation of
$\alpha_{\nu,\mu}$, $\alpha_{\nu,\mu}'$, 
 $\log\left(\tilde{P}_{\nu}^{-\mu}(t)\right) - \nu$ 
and $\log\left(\tilde{Q}_{\nu}^{-\mu}(t)\right) + \nu$ 
with roughly double precision accuracy (see the experiments
of Section~\ref{section:experiments}).
The code was written in Fortran with OpenMP extensions
and compiled with version 4.8.4 of the GNU Fortran compiler.
It was executed on a computer equipped with $28$ Intel Xeon
E5-2697  processor cores running at 2.6 GHz.  The construction
of the table took approximately 24 hours on this machine.

Here, we describe only the construction of the functions $A_1$, $B_1$, $C_1$
and $D_1$.  The construction of the others is extremely similar.
The procedure proceeded in four stages as follows:

\vskip 1em
{\it Stage one: construction of the phase functions and logarithms}

We began this stage of the procedure by constructing a partition 
\begin{equation}
\xi_1 < \xi_2 < \ldots < \xi_{11}
\label{expansions:partition1}
\end{equation}
which divides the interval
\begin{equation}
\left[
\frac{1}{\maxdegree}, \frac{1}{10}
\right]
\label{expansions:interval}
\end{equation}
over which $\frac{1}{\nu}$ is allowed to vary into
$10$ subintervals.  The precise locations of the nodes
$\xi_j$ are not critically important; reasonable
choices were arrived at quickly through trial and error.
Next, we constructed a partition
\begin{equation}
0 = \tau_1 < \tau_2  < \ldots < \tau_{18} = 1
\label{expansions:partition2}
\end{equation}
which divides the interval $[0,1]$ into $17$ subintervals.  Again,
the precise distribution of the nodes $\tau_j$ is not critical
and reasonable choices were arrived at quickly through trial and error.

For each $i=1,\ldots,10$ and $j=1,\ldots,17$, we processed the 
tensor product of  intervals
 $\left[\xi_i,\xi_{i+1}\right] \times\left[\tau_j,\tau_{j+1}\right]$ as follows.
We let $x_1,\ldots,x_{31}$ be the nodes of the $31$-point
Chebyshev grid on the interval $\left[\xi_i,\xi_{i+1}\right]$ and
$y_1,\ldots,y_{31}$ the nodes of 
the $31$-point Chebyshev grid on the interval
$\left[\tau_j,\tau_{j+1}\right]$.  For each pair
$x_k$, $y_l$, the
algorithm of Section~\ref{section:solver} was used
to  calculate $\alpha_{\nu,\mu}$,   $\alpha_{\nu,\mu}'$,
$\log\left(\tilde{P}_\nu^{-\mu}(t)\right) - \nu$, and
$\log\left(\tilde{Q}_\nu^{-\mu}(t)\right) + \nu$
with $\nu$ and $\mu$ taken to be 
\begin{equation}
\nu = \frac{1}{x_k} \ \ \mbox{and} \ \ \mu = 1 + (\nu-1) y_l.
\end{equation}
The precision for the computations was 
$\epsilon = 10^{-17}$.
The functions $\alpha_{\nu,\mu}$ and $\alpha_{\nu,\mu}'$
 were represented
as $30$th order piecewise Chebyshev expansions
on some adaptively determined collection of subintervals of 
$\left[t_{\nu,\mu}^*, \frac{\pi}{2}\right]$,
while the functions
$\log\left(\tilde{P}_\nu^{-\mu}(t)\right) - \nu$ and
$\log\left(\tilde{Q}_\nu^{-\mu}(t)\right) + \nu$ 
 were represented
as $30$th order piecewise Chebyshev expansions
on some adaptively determined collection of subintervals of 
$\left[0,t_{\nu,\mu}^*\right]$.
Using this data,
the nonoscillatory phase function and its derivative
 can be evaluated for any triple $(\nu,\mu,t)$ in the region $\mathcal{O}_1$
via Chebyshev interpolation.
Likewise, the logarithms of the associated Legendre functions
can be evaluated at any point in $\mathcal{N}_1$.

In this stage of the procedure, the differential equation (\ref{introduction:diffeq}) 
was solved
via the algorithm of Section~\ref{section:solver} for
$163,370$ different pairs of the parameters $(\nu,\mu)$,
many of which were large in magnitude.  Obviously, 
this was only possible because our solver runs in 
time independent of $\nu$ and $\mu$.

\vskip 1em
{\it Stage two: formation of unified discretizations}

For each pair of points $\xi$ and $\tau$ such that
$\xi$ is one of the Chebyshev nodes in one
of the subintervals defined by the partition 
(\ref{expansions:partition1})
and $\tau$ is one of the Chebyshev nodes
in one of the subintervals
defined by the partition
(\ref{expansions:partition2}), we used the procedure
of Section~\ref{preliminaries:adaptive} to 
adaptively form discretization of the functions
\begin{equation}
f_{\xi,\tau}(u) = \alpha_{\nu,\mu}'(t(u))
\end{equation}
and
\begin{equation}
g_{\xi,\tau}(u) = \alpha_{\nu,\mu}(t(u)),
\end{equation}
where 
\begin{equation}
\nu = \frac{1}{\xi},  \ \ \mu = 1 + (\nu-1) \tau
\ \ \mbox{and}\ \ 
t(u) = t_{\nu,\mu}^* + \left(\frac{\pi}{2} - t_{\nu,\mu}^* \right)u.
\label{expansions:mapping3}
\end{equation}
The functions $\alpha_{\nu,\mu}$ and $\alpha_{\nu,\mu}'$ 
are evaluated via Chebyshev interpolation using the data constructed in the first stage
of these calculations.  We requested $\epsilon = 10^{-17}$ accuracy
and took the parameter $n$ to be 30.   
For each $\xi$ and $\tau$ considered, this results in a collection
of subintervals of $[0,1]$ on which $f_{\xi,\tau}$ 
is represented with relative accuracy roughly $\epsilon$ via a 
$30$th order Chebyshev expansion
and
another collection of subintervals of $[0,1]$ on which
 $g_{\xi,\tau}$
is represented with relative accuracy roughly $\epsilon$ via a 
$30$th order Chebyshev expansion.  We then formed a unified
discretization
\begin{equation}
\left[a_0,a_1\right],\
\left[a_1,a_2\right],\
\left[a_2,a_3\right],\
\ldots \left[a_{20},a_{21}\right]
\label{expansions:unified1}
\end{equation}
of $[0,1]$ by merging these discretizations; that is, by ensuring
that each subinterval in the discretization of one of the functions
$f_{\xi,\tau}$ or $g_{\xi,\tau}$ is the union of some set of subintervals
of (\ref{expansions:unified1}).

A unified discretization 
\begin{equation}
\left[b_0,b_1\right],\
\left[b_1,b_2\right],\
\left[b_2,b_3\right],\
\ldots \left[b_{17},b_{18}\right]
\label{expansions:unified2}
\end{equation}
for the functions
\begin{equation}
\tilde{f}_{\xi,\tau}(u) = \log\left(\tilde{P}_\nu^{-\mu}(t)\right) - \nu
\end{equation}
and
\begin{equation}
\tilde{g}_{\xi,\tau}(u) = \log\left(\tilde{Q}_\nu^{-\mu}(t)\right) + \nu
\end{equation}
with  $\nu$, $\mu$ and $t$ related to
$\xi$, $\tau$ and $u$ via (\ref{expansions:mapping3})
was formed in the same fashion.

\vskip 1em
{\it Stage three: Construction of the functions $A_1$ and $B_1$}

The function $A_1$ is defined via the formula
\begin{equation}
A_1(\xi,\tau,u) =  \frac{1}{\nu} \alpha_{\nu,\mu}(t),
\end{equation}
where
\begin{equation}
\nu = \frac{1}{\xi}, \ \
\mu = 1 + (\nu-1) \tau\ \ \mbox{and}\ \ 
t = t_{\nu,\mu}^* + \left(\frac{\pi}{2} - t_{\nu,\mu}^* \right)u. 
\label{expansions:mapping}
\end{equation}
Likewise, $B_1$ is defined via
\begin{equation}
B_1(\xi,\tau,u) =  \frac{1}{\nu} \alpha_{\nu,\mu}'(t),
\end{equation}
with $\nu$, $\mu$ and $t$ given by (\ref{expansions:mapping}).  In this
way, we ensure that $A_1$ and $B_1$ are defined on the rectangular
prism
\begin{equation}
\left[\frac{1}{\maxdegree},\frac{1}{10}\right]
\times
\left[0,1\right]
\times
\left[0,1\right],
\label{expansions:prism}
\end{equation}
and hence suitable for representation via a collection of piecewise
trivariate Chebyshev expansions.  

For each $i=1,\ldots,10$, $j=1,\ldots,17$ and 
$k=1,\ldots,21$, we formed the  $30$th order compressed  trivariate 
Chebyshev expansions (as defined in 
Section~\ref{preliminaries:compressed}) for the functions $A_1$ and $B_1$
on the rectangular prism
\begin{equation}
\left[\xi_i,\xi_{i+1}\right] \times \left[\tau_j,\tau_{j+1}\right]
\times \left[a_k,a_{k+1}\right].
\end{equation}
There are $3\sep,570$ such rectangular prisms.
Since an  uncompressed $30$th order trivariate Chebyshev expansion
has  $29\sep,791$ coefficients, 
a  total of 
$212\sep,707\sep,740$ coefficients 
 would be required
to store the uncompressed Chebyshev expansions
of the functions $A_1$ and $B_1$.  
If each coefficient were stored as an IEEE double precision number,
 roughly $1.5$ GB of memory would be required  to store these expansions.
Fortunately, the compressed $30$th order Chebyshev expansions
were far more efficient.  The 
compressed Chebyshev expansions for $A_1$  and $B_1$
had only
$7\sep,839\sep,620$  coefficients.

\vskip 1em
{\it Stage four: construction of the functions $C_1$ and $D_1$}

The function $C_1$ is defined via
\begin{equation}
C_1(\xi,\tau,u) = \log\left(\tilde{P}_\nu^{-\mu}(t) \right) - \nu,
\end{equation}
where
\begin{equation}
\nu = \frac{1}{\xi}, \ \
\mu = 1 + (\nu-1) \tau\ \ \mbox{and}\ \ 
t = \frac{t_{\nu,\mu}^* }{100} + 
\left( t_{\nu,\mu}^*   - \frac{t_{\nu,\mu}^* }{100}\right) u.
\label{expansions:mapping2}
\end{equation}
Finally, $D_1$ is defined via
\begin{equation}
D_1(\xi,\tau,u) = \log\left(\tilde{Q}_\nu^{-\mu}(t) \right) + \nu,
\end{equation}
with $\nu$, $\mu$ and $t$ as in (\ref{expansions:mapping2}).
Obviously, $C_1$ and $D_1$ are also given
on the rectangular prism (\ref{expansions:prism}).

For each $i=1,\ldots,10$, $j=1,\ldots,17$ and 
$k=1,\ldots,17$, we formed the  $30$th order compressed  trivariate 
Chebyshev expansions (as defined in 
Section~\ref{preliminaries:compressed}) for the functions $C_1$ and $D_1$
on the rectangular prism
\begin{equation}
\left[\xi_i,\xi_{i+1}\right] \times \left[\tau_j,\tau_{j+1}\right]
\times \left[b_k,b_{k+1}\right].
\end{equation}
There are $2\sep,890$ such rectangular prisms
and s  total of 
$172\sep,191\sep,980$ coefficients
 would be required
to store the uncompressed Chebyshev expansions
of the functions $C_1$ and $D_1$.  
If each coefficient were stored as an IEEE double precision number,
 roughly $1.3$ GB of memory would be required to store these expansions.
The compressed $30$th order Chebyshev expansions
were far more efficient.  They required only
$4\sep,441\sep,063$
  coefficients 
to store  $C_1$  and $D_1$.

\label{section:expansions}
\end{section}

\begin{section}{An algorithm for the rapid numerical evaluation of the associated Legendre functions}

In this section, we describe the operation of our code for evaluating
the associated Legendre functions $\tilde{P}_\nu^{-\mu}(t)$ and 
$\tilde{Q}_\nu^{-\mu}(t)$ when
\begin{equation}
0 \leq \nu \leq \maxdegree,\   0 \leq \mu \leq \nu
\ \ \mbox{and}\ \  0 < t \leq \frac{\pi}{2}.
\label{algorithm:conditions}
\end{equation}
The code was written in Fortran and its interface to the user consists of two
subroutines, one called {\tt alegendre\_eval_init} and the other {\tt alegendre\_eval}.
The {\tt alegendre\_eval\_init} routine reads the precomputed table constructed
via the procedure of Section~\ref{section:expansions}  from the disk into memory.
The precomputed table used in the  experiments described in this paper
is approximately \tablesize\ megabytes in size.
Once the precomputed table has been loaded, the {\tt alegendre\_eval} can be called.  
It takes as input a triple $(\nu,\mu,t)$ satisfying the conditions
(\ref{algorithm:conditions}).
When $(\nu,\mu,t)$ is in the oscillatory region $\mathcal{O}$, 
it returns the values of $\alpha_{\nu,\mu}(t)$ and $\alpha_{\nu,\mu}'(t)$ as well as
those of 
$\tilde{P}_\nu^{-\mu}$ and  $\tilde{Q}_\nu^{-\mu}$.
  When $(\nu,\mu,t)$ is in the 
nonoscillatory region $\mathcal{N}$, it returns the values
of $\log\left(\tilde{P}_\nu^{-\mu}(t)\right)$ and $\log\left(\tilde{Q}_\nu^{-\mu}(t)\right)$ as well as those
of 
$\tilde{P}_\nu^{-\mu}(t)$ and  $\tilde{Q}_\nu^{-\mu}(t)$.
  Of course, when $t$ is close to $0$,
the latter values might not be representable via the IEEE double format
arithmetic.  In this event, $0$ is returned for
$\tilde{P}_\nu^{-\mu}(t)$ and  $\infty$ for $\tilde{Q}_\nu^{-\mu}(t)$.

The {\tt alegendre\_eval} code is available from the GitHub repository at address

\hskip 2em {\tt http://github.com/JamesCBremerJr/ALegendreEval}.



It  uses several different methods
to evaluate the associated Legendre functions and the associated
auxiliary functions, depending on the values 
of $\nu$, $\mu$  and $t$.  The following description of 
the operation of the {\tt alegendre_eval} code is organized
by listing each such method.

\vskip 1em
{\it Method one: series expansions for $\tilde{P}_\nu^{-\mu}(t)$ and  $\tilde{Q}_\nu^{-\mu}(t)$}

This method is used when
$\nu < 2$ and  $(\nu,\mu,t)$ is in the oscillatory region $\mathcal{O}$.

It consists of evaluating $\tilde{P}_\nu^{-\mu}(t)$ via a   
 truncation of the  series expansion
(\ref{preliminaries:series:2}) 
and  evaluating $\tilde{Q}_\nu^{-\mu}(t)$ via formula (\ref{preliminaries:series:4}).
   As discussed in Section~\ref{preliminaries:series},
when $\mu$  is close to or coincides with an integer, Chebyshev interpolation
in the variable $\mu$ is used to avoid roundoff error in the evaluation
of   (\ref{preliminaries:series:4}).  The value of $\alpha_{\nu,\mu}'(t)$ is calculated 
via  (\ref{preliminaries:alphap}) and $\alpha_{\nu,\mu}$ is computed
using the formula
\begin{equation}
\alpha_{\nu,\mu}
=
\mbox{Arg}\left(\
\tilde{P}_\nu^{-\mu}(t)
+ i \tilde{Q}_\nu^{-\mu}(t)
\right) + 2\pi,
\label{algorithm:1}
\end{equation}
where $\mbox{Arg}(z)$ denotes the principal value of the argument of
the complex number $z$.  The limitation on the range of parameters
for which this method is used
 ensures that the principal branch of the  argument
function is the correct one.

\vskip 1em
{\it Method two: series expansions for $\log\left(\tilde{P}_\nu^{-\mu}(t)\right)$ and  
$\log\left(\tilde{Q}_\nu^{-\mu}(t)\right)$}

This method is used when $\nu < 10$ and $t$  is in the nonoscillatory regime,
and when $10 \leq \nu < 10\sep,000$ and 
\begin{equation}
0 < t < \frac{t_{\nu,\mu}^*}{100}.
\end{equation}

It consists of  evaluating $\log\left(\tilde{P}_\nu^{-\mu}(t)\right)$
via a truncation of  (\ref{preliminaries:series:2.5})  and $\log\left(\tilde{Q}_\nu^{-\mu}(t)\right)$
via  (\ref{preliminaries:series:6}).
 When $\mu$ coincides with or is close to an integer, Chebyshev interpolation
in the parameter $\mu$ is used in the evaluation of 
 (\ref{preliminaries:series:6}).
The values of 
$\tilde{P}_\nu^{-\mu}(t)$
and
$\tilde{Q}_\nu^{-\mu}(t)$
are computed from their logarithms in the obvious fashion.

\vskip 1em
{\it Method three: Macdonald's asymptotic expansions for $\log\left(\tilde{P}_\nu^{-\mu}(t)\right)$ and  
$\log\left(\tilde{Q}_\nu^{-\mu}(t)\right)$}

This method is used when $\nu \geq 10\sep,000$  and
\begin{equation}
0 < t < \frac{t_{\nu,\mu}^*}{100}.
\end{equation}

It consists of  evaluating $\log\left(\tilde{P}_\nu^{-\mu}(t)\right)$
and $\log\left(\tilde{Q}_\nu^{-\mu}\right)$ via Macdonald's asymptotic
expansions (see Section~\ref{preliminaries:macdonald}).
The values of 
$\tilde{P}_\nu^{-\mu}(t)$
and
$\tilde{Q}_\nu^{-\mu}(t)$
are computed from their logarithms in the obvious fashion.

\vskip 1em
{\it Method four:  precomputed expansions}

In all other cases, the precomputed expansions
of the functions
 $A_1, \ldots, A_4$, $B_1,\ldots,B_4$, $C_1,\ldots,C_4$,
whose construction is described
in Section~\ref{section:expansions},
are used to evaluate  $\tilde{P}_\nu^{-\mu}(t)$,
$\tilde{Q}_\nu^{-\mu}(t)$ and the appropriate auxiliary functions.
Here, we describe the use of the functions $A_1$ and $B_1$ 
to evaluate $\alpha_{\nu,\mu}$ and $\alpha_{\nu,\mu}'$ 
in the event that $(\nu,\mu,t)$ is in the 
set $\mathcal{O}_1$. The other cases are extremely similar.

First, we let $\xi = \frac{1}{\nu}$, $\tau = \frac{\mu-1}{\nu-1}$ and
\begin{equation}
u = \frac{t-t^*_{\nu,\mu}}{\frac{\pi}{2}-t^*_{\nu,\mu}}.
\end{equation}
That is, we compute the values of $\xi$, $\tau$ and $u$ defined
by the mapping (\ref{expansions:mapping}) given $\nu$, $\mu$ and $t$.
Next, we  find the smallest positive integer $i$ such that $\xi_i \leq \xi \leq \xi_{i+1}$,
where $\xi_1,\ldots,\xi_{11}$ are the nodes
of the partition (\ref{expansions:partition1}),
the smallest positive integer $j$ such that $\tau_j \leq \tau \leq \tau_{i+1}$,
where $\tau_1,\ldots,\tau_{17}$ are the nodes
of the partition (\ref{expansions:partition2}),
and  the smallest positive integer $k$ such that
$a_k \leq u \leq a_{k+1}$, where
$a_1,a_2,\ldots,a_{21}$ are
the nodes of the partition (\ref{expansions:unified1}). 

Having discovered that $(\xi,\tau,u)$ is in the set
$\left[\xi_i,\xi_{i+1}\right] \times 
\left[\tau_j,\xi_{j+1}\right]\times
\left[a_k,a_{k+1}\right]$, we evaluate the compressed trivariate 
Chebyshev expansions representing $A_1$ and $C_1$ on 
this rectangular prism.  We scale the results
by $\nu$ to obtain the values of $\alpha_{\nu,\mu}$
and $\alpha_{\nu,\mu}'$.  The values of 
$\tilde{P}_\nu^{-\mu}$ and  $\tilde{Q}_\nu^{-\mu}$
are then calculated via 
(\ref{preliminaries:u}) and
(\ref{preliminaries:v}).

\label{section:algorithm}
\end{section}

\begin{section}{Numerical experiments}

In this section, we present the results of numerical experiments which were 
conducted to assess the performance of the {\tt alegendre_eval} routine. 
The task of constructing reference values with which to compare our results
was quite challenging.   All existing packages
of which the author is aware were 
prohibitively slow when evaluating associated Legendre functions
with large noninteger parameters, and existing asymptotic expansions
are either not viable (e.g., the Liouville-Green expansions
(\ref{introduction:BoydDunsterp}) and (\ref{introduction:BoydDunsterq}) 
whose coefficients cannot be readily computed) or only applicable in the case of 
an extremely limited range of parameters (e.g., the trigonometric
expansions  (\ref{introduction:trigp}) and
(\ref{introduction:trigq}) which are catastrophically unstable
even for relatively small values of $\mu$).
As a result, we were quite limited in the extent to which we could verify
our approach in the case of large noninteger parameters.

In the case of integer values of the parameters,
the well-known three term recurrence relations can be used to evaluate
the associated Legendre function accurately, provided extended
precision arithmetic is used to perform the computations.
Consequently, we were able to test our code quite thoroughly in
the case of integer parameters.
We note that the time required to evaluate the associated Legendre functions
using the recurrence relations grows with the magnitudes of the parameters,
making such an approach uncompetitive with the algorithm
of this paper in many cases.

These experiments were carried out on a laptop computer equipped
with an Intel Core i7-5600U processor running at 2.6 GHz and 16 GB of memory.
Our code was compiled with the GNU Fortran compiler version 5.2.1 using
the ``-Ofast'' compiler optimization flag.

\begin{subsection}{The accuracy with which $\alpha_{\nu,\mu}'$ is evaluated for small
noninteger values of $\nu$} 

In these experiments, we measured the accuracy with which {\tt alegendre_eval}
calculates  $\alpha_{\nu,\mu}'$ in the oscillatory region.
Reference values were calculated using version 11 of 
Wolfram's  Mathematica package.  The cost of the reference calculations
was prohibitively expensive for large $\nu$, with the consequence
that we only considered values of $\nu$ between $0$ and $1\sep,000$.

In each experiment, we choose $10$ pairs $(\nu,\mu)$ by
first picking a random value of $\nu$ in a given range,
and then choosing  a random value of $\mu$ in the interval $(0,\nu)$.
For each pair chosen in this fashion, we evaluated 
$\alpha_{\nu,\mu}'$ at $100$ equispaced points
 either in  interval
\begin{equation}
\left(\arcsin\left(\frac{\sqrt{\mu^2-\frac{1}{4}}}{\nu+\frac{1}{2}}\right),\frac{\pi}{2}\right)
\end{equation}
or in the interval
\begin{equation}
\left(\frac{1}{1000},\frac{\pi}{2}\right),
\end{equation}
depending on whether $\mu > \frac{1}{2}$ or not.
Table~\ref{table1} reports the results.  
There, each row corresponds
to one experiment and gives 
the largest relative error observed in $\alpha_{\nu,\mu}'$ 
as well as the average time taken by the {\tt alegendre_eval} routine.

\begin{table}[h!]
\begin{center}
\small
 \begin{tabular}{lcc}
 \toprule
 Range of $\nu$ & Maximum relative               & Average evaluation \\
                & error in $\alpha_{\nu,\mu}'$   & time (in seconds)   \\
 \midrule
0 - 1 & 2.26\e{-14} & 2.89\e{-06}  \\
1 - 5 & 2.62\e{-15} & 1.92\e{-06}  \\
5 - 10 & 2.38\e{-15} & 1.65\e{-06}  \\
10 - 50 & 4.15\e{-15} & 3.13\e{-06}  \\
50 - 100 & 8.53\e{-15} & 2.15\e{-06}  \\
100 - 500 & 1.88\e{-14} & 2.46\e{-06}  \\
500 - 1\sep,000 & 3.49\e{-14} & 1.45\e{-06}  \\
 \bottomrule
 \end{tabular}

\end{center}
\vskip .5em
\caption{The results of the experiments of Section~\ref{experiments:1} in which
the accuracy with which {\tt alegendre_eval} calculates
$\alpha_{\nu,\mu}'$ for small values of $\nu$  is tested via comparison
with Wolfram's Mathematica package.
}
\label{table1}
\end{table}

\label{experiments:1}
\end{subsection}


\begin{subsection}{The accuracy with which logarithms are evaluated in the case of 
small noninteger values of $\nu$} 

In these experiments, we measured the accuracy with which {\tt alegendre_eval}
calculates  the functions
\begin{equation}
\log\left(\tilde{P}_\nu^{-\mu}(t)\right) - \nu
\ \ \mbox{and} \ \ 
\log\left(\tilde{Q}_\nu^{-\mu}(t)\right) + \nu
\label{experiments:1000}
\end{equation}
in the nonoscillatory region.
High accuracy reference values for these experiments were calculated using version 11 of 
Wolfram's  Mathematica package.  Again the high cost of the reference calculations
led us to only consider values of $\nu$ between $0$ and $1\sep,000$.

In each experiment, we choose $10$ pairs $(\nu,\mu)$ by
first picking a random value of $\nu$ in a given range,
and then choosing  a random value of $\mu$ in the interval $\left(\frac{1}{2},\nu\right)$.
For each pair chosen in this fashion, we evaluated 
the functions (\ref{experiments:1000})
 at $100$ equispaced points in the interval
\begin{equation}
\left(0,\arcsin\left(\frac{\sqrt{\mu^2-\frac{1}{4}}}{\nu+\frac{1}{2}}\right)\right).
\end{equation}
Table~\ref{table2} reports the results.  
There, each row corresponds to one experiment and gives 
the largest relative error observed in each of the functions
(\ref{experiments:1000}),
as well as the average time taken by the {\tt alegendre_eval} routine.

\begin{table}[h!]
\begin{center}
\small
 \begin{tabular}{lccc}
 \toprule
 Range of $\nu$ & Maximum relative & Maximum relative & Average evaluation \\
                & error in         & error in         & time (in seconds)   \\
 &$\log(\tilde{P}_\nu^{-\mu}(t))-\nu$ &$\log(\tilde{Q}_\nu^{-\mu}(t))+\nu$ & \\
 \midrule
0.5
 - 1 & 3.36\e{-16} & 2.58\e{-15} & 1.22\e{-06}  \\
1 - 5 & 3.21\e{-16} & 9.28\e{-16} & 1.32\e{-06}  \\
5 - 10 & 8.85\e{-16} & 9.14\e{-15} & 1.59\e{-06}  \\
10 - 50 & 4.39\e{-15} & 4.43\e{-15} & 3.34\e{-06}  \\
50 - 100 & 2.58\e{-15} & 3.49\e{-15} & 1.61\e{-06}  \\
100 - 500 & 4.21\e{-15} & 4.47\e{-15} & 2.56\e{-06}  \\
500 - 1\sep,000 & 2.54\e{-15} & 3.24\e{-15} & 1.70\e{-06}  \\
 \bottomrule
 \end{tabular}

\end{center}
\vskip .5em
\caption{The results of the experiments of Section~\ref{experiments:2} in which
the accuracy with which {\tt alegendre_eval} calculates
$\alpha_{\nu,\mu}'$ for small values of $\nu$  is tested via comparison
with Wolfram's Mathematica package.}
\label{table2}
\end{table}

\label{experiments:2}
\end{subsection}


\begin{subsection}{The accuracy with which $\alpha_{\nu,\mu}'$ is evaluated
in the case of large $\nu$ and small $\mu$}

In these experiments, we measured the accuracy with which {\tt alegendre_eval}
calculates  $\alpha_{\nu,\mu}'$ in the oscillatory region by comparison
with values obtained using the trigonometric expansions (\ref{introduction:trigp}) 
and (\ref{introduction:trigq}).
Since these expansions are numerically unstable, and catastrophically so when $\mu$ is large,
we considered only pairs of the parameters $(\nu,\mu)$ with $\mu$ small in magnitude.
  Even so, $1\sep,000$ digit arithmetic was required in
order to obtain accurate reference values for these experiments.

In each experiment, we choose $10$ pairs $(\nu,\mu)$ by
first picking a random value of $\nu$ in a given range, and then 
choosing  a random value of $\mu$ in the interval $\left(\frac{1}{2},\frac{\nu}{100}\right)$.
For each pair, the function $\alpha_{\nu,\mu}'$ was evaluated
at $100$ equispaced points in the interval
\begin{equation}
\left(
\max\left(
\arcsin\left(
\frac{\sqrt{\mu^2-\frac{1}{4}}}{\nu+\frac{1}{2}}\right),\frac{\pi}{6}\right),
\frac{\pi}{2}\right).
\label{experiments:1001}
\end{equation}
We note that the trigonometric expansions used here only converge in the
interval $\left(\frac{\pi}{6},\frac{5\pi}{6}\right)$, hence
the choice of the interval (\ref{experiments:1001}).

\begin{table}[h!]
\begin{center}
\small
 \begin{tabular}{lcc}
 \toprule
 Range of $\nu$ & Maximum relative               & Average evaluation \\
                & error in $\alpha_{\nu,\mu}'$   & time (in seconds)   \\
 \midrule
1\sep,000 - 5\sep,000 & 1.91\e{-15} & 1.53\e{-06}  \\
5\sep,000 - 10\sep,000 & 1.41\e{-15} & 1.00\e{-06}  \\
10\sep,000 - 50\sep,000 & 1.05\e{-15} & 1.06\e{-06}  \\
50\sep,000 - 100\sep,000 & 8.69\e{-16} & 7.72\e{-07}  \\
100\sep,000 - 500\sep,000 & 7.30\e{-16} & 6.91\e{-07}  \\
500\sep,000 - 1\sep,000\sep,000 & 8.15\e{-16} & 5.46\e{-07}  \\
 \bottomrule
 \end{tabular}

\end{center}
\vskip .5em
\caption{The results of the experiments of Section~\ref{experiments:3} in which
the accuracy with which {\tt alegendre_eval} calculates
$\alpha_{\nu,\mu}'$ for large $\nu$ and small $\mu$ is tested via comparison
with the trigonometric expansions (\ref{introduction:trigp})
and (\ref{introduction:trigq}).
}
\label{table3}
\end{table}

\label{experiments:3}
\end{subsection}


\begin{subsection}{The accuracy with which $\alpha_{\nu,\mu}'$ 
is evaluated in the case  of integer parameters}

In these experiments, the accuracy with which $\alpha_{\nu,\mu}'$
is evaluated in the oscillatory regime was measured by comparison
with reference values calculated using the well-known three
term recurrence relations satisfied by the associated Legendre
functions.     The reference calculations were conducted in
extended precision arithmetic in order to ensure accuracy.

The experiments of this section proceeded just as those described in Section~\ref{experiments:1},
except only integer values of the parameters were considered.
Table~\ref{table4} displays the results.

\begin{table}[h!]
\begin{center}
\small
 \begin{tabular}{lcc}
 \toprule
 Range of $\nu$ & Maximum relative               & Average evaluation \\
                & error in $\alpha_{\nu,\mu}'$   & time (in seconds)   \\
 \midrule
10 - 50 & 2.35\e{-14} & 2.79\e{-06}  \\
50 - 100 & 4.71\e{-15} & 1.87\e{-06}  \\
100 - 500 & 4.96\e{-15} & 2.86\e{-06}  \\
500 - 1\sep,000 & 2.86\e{-14} & 1.52\e{-06}  \\
1\sep,000 - 5\sep,000 & 8.62\e{-15} & 1.37\e{-06}  \\
5\sep,000 - 10\sep,000 & 5.94\e{-15} & 1.14\e{-06}  \\
10\sep,000 - 50\sep,000 & 2.74\e{-14} & 1.25\e{-06}  \\
50\sep,000 - 100\sep,000 & 7.36\e{-14} & 9.81\e{-07}  \\
100\sep,000 - 500\sep,000 & 1.86\e{-14} & 9.71\e{-07}  \\
500\sep,000 - 1\sep,000\sep,000 & 3.09\e{-14} & 8.86\e{-07}  \\
 \bottomrule
 \end{tabular}

\end{center}
\vskip .5em
\caption{The results of the experiments of Section~\ref{experiments:4} in which
the accuracy with which {\tt alegendre_eval} calculates
$\alpha_{\nu,\mu}'$ for integers values of the parameters is tested.
}
\label{table4}
\end{table}

\label{experiments:4}
\end{subsection}


\begin{subsection}{The accuracy with which the logarithms
are evaluated in the case  of integer parameters}

In these experiments, we measured the accuracy with which
{\tt alegendre_eval} calculates
the functions (\ref{experiments:1000}) in the nonoscillatory
regime.
Reference values were calculated using the well-known three
term recurrence relations satisfied by the associated Legendre
functions.   The reference calculations were conducted in
extended precision arithmetic in order to ensure accuracy.

These experiments proceeded just as those described in Section~\ref{experiments:2},
except only integer values of the parameters were considered.
Table~\ref{table5} displays the results.

\begin{table}[h!]
\begin{center}
\small
 \begin{tabular}{lccc}
 \toprule
 Range of $\nu$ & Maximum relative & Maximum relative & Average evaluation \\
                & error in         & error in         & time (in seconds)   \\
 &$\log(\tilde{P}_\nu^{-\mu}(t))-\nu$ &$\log(\tilde{Q}_\nu^{-\mu}(t))+\nu$ & \\
 \midrule
10 - 50 & 4.21\e{-15} & 4.65\e{-15} & 3.59\e{-06}  \\
50 - 100 & 3.42\e{-15} & 3.32\e{-15} & 2.22\e{-06}  \\
100 - 500 & 3.07\e{-15} & 4.07\e{-15} & 2.28\e{-06}  \\
500 - 1\sep,000 & 2.95\e{-15} & 3.01\e{-15} & 1.88\e{-06}  \\
1\sep,000 - 5\sep,000 & 2.63\e{-15} & 4.14\e{-15} & 1.55\e{-06}  \\
5\sep,000 - 10\sep,000 & 1.98\e{-15} & 1.83\e{-15} & 9.60\e{-07}  \\
10\sep,000 - 50\sep,000 & 1.98\e{-15} & 2.68\e{-15} & 1.69\e{-06}  \\
50\sep,000 - 100\sep,000 & 1.63\e{-15} & 2.07\e{-15} & 1.17\e{-06}  \\
100\sep,000 - 500\sep,000 & 1.73\e{-15} & 1.63\e{-15} & 1.21\e{-06}  \\
500\sep,000 - 1\sep,000\sep,000 & 1.67\e{-15} & 2.23\e{-15} & 1.16\e{-06}  \\
 \bottomrule
 \end{tabular}

\end{center}
\vskip .5em
\caption{The results of the experiments of Section~\ref{experiments:5} in which
the accuracy with which {\tt alegendre_eval} calculates
the functions $\log\left(\tilde{P}_\nu^{-\mu}(t)\right)-\nu$ 
and $\log\left(\tilde{Q}_\nu^{-\mu}(t)\right)-\nu$
for integers values of the parameters is tested.
}
\label{table5}
\end{table}

\label{experiments:5}
\end{subsection}


\begin{subsection}{The accuracy with which the associated Legendre
functions are evaluated in the case  of integer parameters}

In these experiments, we measured the accuracy with which
{\tt alegendre_eval} calculates the functions $\tilde{P}_\nu^{-\mu}$
and $\tilde{Q}_\nu^{-\mu}$ 
 in the case of integer values of the parameters.
Reference values were calculated using the  three term recurrence relations.
As usual, extended precision arithmetic was used during the reference
calculations in order to ensure their accuracy.

In each experiment, $10$ pairs of the parameters $(\nu,\mu)$ were constructed 
by first choosing an integer value of $\nu$ in a given range at random
and then choosing  an integer value of $\mu$ in the range $(0,\nu)$
at random.   For each such pair, we evaluated  the function
$\tilde{P}_\nu^{-\mu}(t) + i\tilde{Q}_{\nu}^{-\mu}(t)$ at $100$ equispaced
points either in  interval
\begin{equation}
\left(\arcsin\left(\frac{\sqrt{\mu^2-\frac{1}{4}}}{\nu+\frac{1}{2}}\right),\frac{\pi}{2}\right)
\end{equation}
or in the interval
\begin{equation}
\left(\frac{1}{1000},\frac{\pi}{2}\right),
\end{equation}
depending on whether $\mu > \frac{1}{2}$ or not.  Table~\ref{table6}
reports the results.  Each row corresponds to one experiment and 
reports the largest relative error which was observed as well as 
the average evaluation time.  We note that we considered
 the function  $\tilde{P}_\nu^{-\mu}(t) + i\tilde{Q}_{\nu}^{-\mu}(t)$ 
 because, unlike 
$\tilde{P}_\nu^{-\mu}(t)$ and $\tilde{Q}_\nu^{-\mu}(t)$,
 its absolute value is nonoscillatory and  does not have roots on the interval
$(0,\pi)$.

\begin{table}[h!]
\begin{center}
\small
 \begin{tabular}{lcc}
 \toprule
 Range of $\nu$ & Maximum relative  & Average evaluation \\
                & error in           & time (in seconds)   \\
 &$\tilde{P}_\nu^{-\mu}(t) + i \tilde{Q}_\nu^{-\mu}(t)$  \\
 \midrule
10 - 50 & 2.62\e{-13} & 4.05\e{-06}  \\
50 - 100 & 4.20\e{-13} & 1.95\e{-06}  \\
100 - 500 & 1.20\e{-12} & 2.24\e{-06}  \\
500 - 1\sep,000 & 1.72\e{-12} & 1.48\e{-06}  \\
1\sep,000 - 5\sep,000 & 8.57\e{-12} & 2.55\e{-06}  \\
5\sep,000 - 10\sep,000 & 1.38\e{-11} & 1.14\e{-06}  \\
10\sep,000 - 50\sep,000 & 8.51\e{-11} & 2.51\e{-06}  \\
50\sep,000 - 100\sep,000 & 9.07\e{-11} & 1.95\e{-06}  \\
100\sep,000 - 500\sep,000 & 9.83\e{-10} & 1.31\e{-06}  \\
500\sep,000 - 1\sep,000\sep,000 & 8.25\e{-10} & 1.20\e{-06}  \\
 \bottomrule
 \end{tabular}

\end{center}
\vskip .5em
\caption{The results of the experiments of Section~\ref{experiments:6} in which
the accuracy with which {\tt alegendre_eval} evaluates
the function $\tilde{P}_\nu^{-\mu}(t) + i \tilde{Q}_\nu^{-\mu}(t)$
for integers values of the parameters is tested.
}
\label{table6}
\end{table}

From Table~\ref{table6}, we see that the relative errors in
the calculated values  of the associated Legendre functions increase
as a function of the parameter $\nu$.  This is expected,
and consistent with the condition number of the evaluation
of the function $\tilde{P}_\nu^{-\mu}(t) + i \tilde{Q}_\nu^{-\mu}(t)$.



\label{experiments:6}
\end{subsection}

\label{section:experiments}
\end{section}

\begin{section}{Conclusions}


Nonoscillatory phase functions provide powerful
theoretical tools for analyzing the solutions of second order
differential equations as well as a framework for the design of simple and efficient
numerical algorithms.    
Here, we have designed a  scheme for the numerical
evaluation of the associated Legendre functions on the cut using
this  framework.  Our approach is simple-minded and highly effective.
Moreover, by  making use of  the algorithms of \cite{BremerKummer} and \cite{Bremer-Rokhlin},
it can be applied in the case of many other special functions
satisfying second order differential equations, such as the 
prolate spheroidal wave functions and the generalized Laguerre functions.
The author will report on the use of the techniques of
this paper to evaluate other special functions of interest
at a later date.

In the nonoscillatory region, our algorithm calculates the logarithms  
of the associated Legendre functions
as well as their values.  This is useful in cases in which the magnitudes
of those functions are too large or too small to be encoded
using the IEEE double precision format.
In the oscillatory region, in addition to the values
of the associated Legendre functions, our algorithm also returns the values
of a nonoscillatory phase function for the associated Legendre
differential equation and its derivative.
This is extremely helpful when  computing the zeros of special functions
\cite{BremerZeros}, and when applying special function transforms
via the butterfly algorithm (see, for instance, 
\cite{Candes-Demanet-Ying1,butterfly2,butterfly1,Candes-Demanet-Ying2,Michielssen,ONeil-Rokhlin}).
The author will report  on the use nonoscillatory phase functions to rapidly compute
the roots of the associated Legendre functions and to rapidly
apply the spherical harmonic transform at a later date.

\label{section:conclusion}
\end{section}

\begin{section}{Acknowledgments}
The author is grateful to Zydrunas Gimbutas of NIST Boulder
for providing his arbitrary precision arithmetic codes
for evaluating the associated Legendre functions
of large integer degrees and orders via the three term recurrence relations
they satisfy.
This work was supported in part by a UC Davis Chancellor's Fellowship.
\end{section}

\bibliographystyle{acm}
\bibliography{alegendre}

\end{document}